\documentclass[11pt]{amsart}
\usepackage{amsmath,amsfonts,amssymb}
\begin{document}

\title[Operator error estimates for homogenization]{Operator error estimates for homogenization
of the elliptic Dirichlet problem in a bounded domain}

\author{M.~A.~Pakhnin}

\keywords{Periodic differential operators, homogenization, effective operator, operator error estimates}

\subjclass[2000]{Primary 35B27}

\address{St. Petersburg State University, Department of Physics, Ul'yanovskaya 3, Petrodvorets,
St.~Petersburg, 198504, Russia}

\email{mpakhnin@yandex.ru}

\author{T.~A.~Suslina}
\thanks{Supported by RFBR (grant no. 11-01-00458-a)
and the Program of support of the leading scientific schools}

\address{St. Petersburg State University, Department of Physics, Ul'yanovskaya 3, Petrodvorets,
St.~Petersburg, 198504, Russia}

\email{suslina@list.ru}

\date{}

\begin{abstract}
Let $\mathcal{O} \subset \mathbb{R}^d$ be a bounded domain of class $C^2$.
In the Hilbert space $L_2(\mathcal{O};\mathbb{C}^n)$, we consider a matrix
elliptic second order differential operator $\mathcal{A}_{D,\varepsilon}$ with the Dirichlet
boundary condition. Here $\varepsilon>0$ is the small parameter.
The coefficients of the operator are periodic and depend on $\mathbf{x}/\varepsilon$.
We find approximation of the operator $\mathcal{A}_{D,\varepsilon}^{-1}$ in the norm of operators
acting from $L_2(\mathcal{O};\mathbb{C}^n)$ to the Sobolev space $H^1(\mathcal{O};\mathbb{C}^n)$
with an error term $O(\sqrt{\varepsilon})$.
This approximation is given by the sum of the operator $(\mathcal{A}^0_D)^{-1}$ and the first order
corrector, where $\mathcal{A}^0_D$ is the effective operator with constant coefficients and
with the Dirichlet boundary condition.
\end{abstract}

\maketitle

\section*{Introduction}

The paper concerns homogenization theory of periodic differential operators (DO's).
A broad literature is devoted to homogenization problems in the small period limit.
First of all, we mention the books [BeLP], [BaPa], [ZhKO].

\noindent\textbf{0.1. Operator-theoretic approach to homogenization problems.}
In a series of papers [BSu1-5] by M.~Sh.~Birman and T.~A.~Suslina
a new operator-theoretic (spectral) approach to homogenization problems was suggested and developed.
By this approach, the so-called operator error estimates in homogenization
problems for elliptic DO's were obtained.
Matrix elliptic DO's acting in $L_2(\mathbb{R}^d;\mathbb{C}^n)$ and admitting a factorization of the form
$\mathcal{A}_\varepsilon = b(\mathbf{D})^* g(\mathbf{x}/\varepsilon)b(\mathbf{D})$, $\varepsilon>0$,
were studied. Here $g(\mathbf{x})$ is an $(m\times m)$-matrix-valued function;
it is assumed to be bounded, uniformly positive definite and periodic with respect to some lattice $\Gamma$.
By $\Omega$ we denote the elementary cell of the lattice $\Gamma$.
It is assumed that $m \geq n$ and $b(\mathbf{D})$ is an $(m\times n)$-matrix homogeneous
first order DO such that $\text{rank}\, b(\boldsymbol{\xi})=n$ for $0 \ne
\boldsymbol{\xi} \in \mathbb{R}^d$.
The simplest example of such operator is the scalar elliptic operator
${\mathcal{A}}_\varepsilon= - \text{div}\, g(\mathbf{x}/\varepsilon) \nabla$.
The operator of elasticity theory also can be written in the required form.
These and other examples are considered in [BSu2] in detail.

In [BSu1-5], the equation
${\mathcal{A}}_\varepsilon \mathbf{u}_\varepsilon + \mathbf{u}_\varepsilon = \mathbf{F}$,
where $\mathbf{F} \in L_2(\mathbb{R}^d;\mathbb{C}^n)$, was considered.
The behavior of the solution $\mathbf{u}_\varepsilon$ for small $\varepsilon$ was studied.
The solution $\mathbf{u}_\varepsilon$ converges in $L_2(\mathbb{R}^d;\mathbb{C}^n)$
to the solution $\mathbf{u}_0$ of the "homogenized"\ equation
${\mathcal{A}}^0 \mathbf{u}_0 + \mathbf{u}_0 = \mathbf{F}$, as $\varepsilon \to 0$.
Here ${\mathcal{A}}^0 = b(\mathbf{D})^* g^0 b(\mathbf{D})$ is the \textit{effective operator}
with the constant effective matrix $g^0$.
In [BSu1,2], it was proved that
$$
\| \mathbf{u}_\varepsilon - \mathbf{u}_0 \|_{L_2(\mathbb{R}^d)}
\leq C \varepsilon \| \mathbf{F} \|_{L_2(\mathbb{R}^d)}.
$$
In operator terms it means that the resolvent
$({\mathcal{A}}_\varepsilon +I)^{-1}$ converges in the operator norm in $L_2(\mathbb{R}^d;\mathbb{C}^n)$
to the resolvent of the effective operator, as $\varepsilon \to 0$, and
$$
\| ({\mathcal{A}}_\varepsilon +I)^{-1} - ({\mathcal{A}}^0 +I)^{-1} \|_{L_2(\mathbb{R}^d) \to L_2(\mathbb{R}^d)}
\le C \varepsilon.
\eqno(0.1)
$$

In [BSu3,4], more accurate approximation of the resolvent
$({\mathcal{A}}_\varepsilon +I)^{-1}$ in the operator norm in $L_2(\mathbb{R}^d;\mathbb{C}^n)$
with an error term $O(\varepsilon^2)$ was obtained.

In [BSu5], approximation of the resolvent $({\mathcal{A}}_\varepsilon +I)^{-1}$
in the norm of operators acting from $L_2(\mathbb{R}^d;\mathbb{C}^n)$ to the Sobolev space
$H^1(\mathbb{R}^d;\mathbb{C}^n)$ was found:
$$
\| ({\mathcal{A}}_\varepsilon +I)^{-1} - ({\mathcal{A}}^0 +I)^{-1} -
\varepsilon K(\varepsilon) \|_{L_2(\mathbb{R}^d) \to H^1(\mathbb{R}^d)} \leq C \varepsilon;
\eqno(0.2)
$$
this corresponds to approximation of $\mathbf{u}_\varepsilon$ in the "energy"\ norm.
Here $K(\varepsilon)$ is a corrector. It contains rapidly oscillating factors and so depends on
$\varepsilon$.

Estimates (0.1), (0.2) are called the \textit{operator error estimates}. They are order-sharp;
the constants in estimates are controlled explicitly in terms of the problem data.
The method of [BSu1--5] is based on the scaling transformation and the Floquet-Bloch theory.
The operator $\mathcal{A} = b(\mathbf{D})^* g(\mathbf{x})b(\mathbf{D})$ is decomposed in the direct
integral of the operators $\mathcal{A}(\mathbf{k})$ acting in $L_2(\Omega;\mathbb{C}^n)$ and
depending on the parameter $\mathbf{k}$ (the quasimomentum).
The operator family $\mathcal{A}(\mathbf{k})$
has discrete spectrum and depends on $\mathbf{k}$ analytically.
It is studied by methods of the analytic perturbation theory.
It turns out that only the spectral characteristics of the operator $\mathcal{A}$
near the bottom of its spectrum are important for constructing the effective operator
and obtaining error estimates. This shows that homogenization can be studied as
a spectral threshold effect.

\smallskip\noindent\textbf{0.2. A different approach} to operator error estimates in homogenization
problems was suggested by V.~V.~Zhikov. In [Zh1, Zh2, ZhPas, Pas], the scalar elliptic operator
$- \text{div}\, g(\mathbf{x}/\varepsilon) \nabla$ (where $g(\mathbf{x})$ is a matrix with
real entries) and the system of elasticity theory were studied.
Estimates of the form (0.1), (0.2) for the corresponding problems in $\mathbb{R}^d$ were obtained.
The method was based on analysis of the first order approximation to the solution and introducing
of an additional parameter (the shift by the vector $\boldsymbol{\omega}\in \Omega$).
Besides the problems in $\mathbb{R}^d$, homogenization problems in
a bounded domain $\mathcal{O} \subset \mathbb{R}^d$ with the Dirichlet or Neumann boundary condition
were studied. Approximation of the solution in $H^1(\mathcal{O})$ was deduced from
the corresponding result in $\mathbb{R}^d$. Due to the "boundary layer"\ influence,
estimates in a bounded domain become worse and the error term is $O(\varepsilon^{1/2})$.
The estimate $\|\mathbf{u}_\varepsilon - \mathbf{u}_0\|_{L_2(\mathcal{O})}
\leq C \varepsilon^{1/2}\|\mathbf{F}\|_{L_2(\mathcal{O})}$
follows from approximation of the solution in $H^1(\mathcal{O})$ by roughening.

Similar results for the operator $- \text{div}\, g(\mathbf{x}/\varepsilon) \nabla$ in a bounded domain
with the Dirichlet or Neumann boundary condition were obtained
in the papers [Gr1, Gr2] by G.~Griso by the "unfolding"\ method.

\smallskip\noindent\textbf{0.3. Main results.}
In the present paper, we study matrix DO's ДО $\mathcal{A}_{D,\varepsilon}$
in a bounded domain $\mathcal{O}\subset \mathbb{R}^d$ of class $C^2$.
The operator $\mathcal{A}_{D,\varepsilon}$ is defined by
the differential expression $b(\mathbf{D})^* g(\mathbf{x}/\varepsilon) b(\mathbf{D})$
with the Dirichlet condition on $\partial \mathcal{O}$.
The effective operator $\mathcal{A}_D^0$ is given by the expression
$b(\mathbf{D})^* g^0 b(\mathbf{D})$ with the Dirichlet boundary condition.
The behavior for small $\varepsilon$ of the solution $\mathbf{u}_\varepsilon$ of the equation
$\mathcal{A}_{D,\varepsilon} \mathbf{u}_\varepsilon =\mathbf{F}$, where
$\mathbf{F} \in L_2(\mathcal{O};\mathbb{C}^n)$, is studied.
Estimates for the $H^1$-norm of the difference of the solution $\mathbf{u}_\varepsilon$
and its first order approximation are obtained.
Roughening this result, we estimate $\|\mathbf{u}_\varepsilon - \mathbf{u}_0\|_{L_2(\mathcal{O})}$.
Here $\mathbf{u}_0$ is the solution of the equation $\mathcal{A}^0_D \mathbf{u}_0 =\mathbf{F}$.

Main results of the paper are Theorems 6.1 and 7.1.
In operator terms, the following estimates are obtained:
$$
\| \mathcal{A}_{D,\varepsilon}^{-1} - (\mathcal{A}^0_D)^{-1} - \varepsilon K_D(\varepsilon) \|_{L_2(\mathcal{O}) \to H^1(\mathcal{O})}
\leq C \varepsilon^{1/2},
\eqno(0.3)
$$
$$
\| \mathcal{A}_{D,\varepsilon}^{-1} - (\mathcal{A}^0_D)^{-1} \|_{L_2(\mathcal{O}) \to L_2(\mathcal{O})}
\leq C \varepsilon^{1/2}.
\eqno(0.4)
$$
Here $K_D(\varepsilon)$ is the corresponding corrector.
The form of the corrector depends on the properties of the
periodic solution $\Lambda(\mathbf{x})$ of the auxiliary problem (1.5).
If $\Lambda$ is bounded, the corrector has a standard form
(Theorem 6.1). In the general case the corrector contains an auxiliary smoothing operator (Theorem 7.1).
Besides approximation of the solution $\mathbf{u}_\varepsilon$ in $H^1(\mathcal{O};\mathbb{C}^n)$,
we also obtain approximation of the "flux" $\mathbf{p}_\varepsilon :=
g^\varepsilon b(\mathbf{D})\mathbf{u}_\varepsilon$ in $L_2(\mathcal{O};\mathbb{C}^m)$.

\smallskip\noindent\textbf{0.4. The method} is based on using estimates (0.1), (0.2)
for homogenization problem in $\mathbb{R}^d$ obtained in [BSu2,5]
and on the tricks suggested in [Zh2], [ZhPas] that allow one to deduce estimate (0.3) from (0.1), (0.2).
Main difficulties are related to estimating of the "discrepancy"\ $\mathbf{w}_\varepsilon$ which satisfies
the equation $\mathcal{A}_\varepsilon \mathbf{w}_\varepsilon =0$ in $\mathcal{O}$ and the boundary
condition $\mathbf{w}_\varepsilon = \varepsilon K_D(\varepsilon)\mathbf{F}$ on $\partial \mathcal{O}$.
Note that we can not use the facts specific for scalar elliptic equations, because we study
a wide class of matrix elliptic DO's.

\smallskip\noindent\textbf{0.5. Error estimates in} $L_2(\mathcal{O})$.
It must be mentioned that estimate (0.4) is quite a rough consequence of (0.3).
So, the refinement of estimate (0.4) is a natural problem.
In [ZhPas], for the case of the scalar elliptic operator
$- \text{div}\, g(\mathbf{x}/\varepsilon) \nabla$ (where $g(\mathbf{x})$ is a matrix with real entries)
an estimate for
$\|\mathcal{A}_{D,\varepsilon}^{-1} - (\mathcal{A}^0_D)^{-1} \|_{L_2\to L_2}$
of order $\varepsilon^{\frac{d}{2d-2}}$ for $d \geq 3$ and of order $\varepsilon |\log \varepsilon|$
for $d=2$ was obtained.
The proof essentially relies on the maximum principle which is specific for scalar
elliptic equations.

Using the results and technique of the present paper, one of the authors has obtained a
\textit{sharp order operator error estimate}
$$
\| \mathcal{A}_{D,\varepsilon}^{-1} - (\mathcal{A}^0_D)^{-1} \|_{L_2(\mathcal{O}) \to L_2(\mathcal{O})}
\leq C \varepsilon.
$$
A separate paper [Su] is devoted to the proof of this result.

\smallskip\noindent\textbf{0.6. The plan of the paper.}
The paper contains seven sections.
In Section 1, the class of operators acting in $L_2(\mathbb{R}^d;\mathbb{C}^n)$
is introduced, the effective operator and the corrector are described, and the needed results
from [BSu2,5] are formulated.
In Section 2, properties of the matrix-valued function
$\Lambda$ are described. In Section 3, we introduce the operator smoothing in Steklov's sense
and prove one more theorem for homogenization problem in $\mathbb{R}^d$.
Section 4 contains the statement of the problem in a bounded domain
and description of the "homogenized"\ problem.
In Section 5, we prove some auxiliary statements needed for further investigation.
Main results of the paper are formulated and proved in Sections 6 and 7.
Herewith, in Section 6 the case where $\Lambda \in L_\infty$ is studied, while
in Section 7 the general case is considered.

\smallskip\noindent\textbf{0.6. Notation.} Let $\mathfrak{H}$ and $\mathfrak{H}_*$ be complex
separable Hilbert spaces. The symbols $(\cdot,\cdot)_{\mathfrak{H}}$ and $\|\cdot\|_{\mathfrak{H}}$
stand for the inner product and the norm in $\mathfrak{H}$;
the symbol $\|\cdot\|_{\mathfrak{H} \to \mathfrak{H}_*}$
denotes the norm of a linear continuous operator acting from $\mathfrak{H}$ to $\mathfrak{H}_*$.

The symbols $\langle \cdot, \cdot \rangle$ and $|\cdot|$ stand for the inner product
and the norm in $\mathbb{C}^n$; $\mathbf{1} = \mathbf{1}_n$ is the identity $(n\times n)$-matrix.
If $a$ is an $(n\times n)$-matrix, the symbol $|a|$ denotes the norm
of the matrix $a$ viewed as a linear operator in $\mathbb{C}^n$.
We use the notation $\mathbf{x} = (x_1,\dots,x_d)\in \mathbb{R}^d$, $iD_j = \partial_j = \partial/\partial x_j$,
$j=1,\dots,d$, $\mathbf{D} = -i \nabla = (D_1,\dots,D_d)$.
The $L_p$-classes of $\mathbb{C}^n$-valued functions in a domain ${\mathcal O} \subset \mathbb{R}^d$
are denoted by $L_p({\mathcal O};\mathbb{C}^n)$, $1 \leq p \leq \infty$.
The Sobolev classes of $\mathbb{C}^n$-valued functions in a domain ${\mathcal O} \subset \mathbb{R}^d$
are denoted by $H^s({\mathcal O};\mathbb{C}^n)$.
By $H^1_0(\mathcal{O};\mathbb{C}^n)$ we denote the closure of $C_0^\infty(\mathcal{O};\mathbb{C}^n)$
in $H^1(\mathcal{O};\mathbb{C}^n)$.
If $n=1$, we write simply $L_p({\mathcal O})$, $H^s({\mathcal O})$, etc., but sometimes
we use such abbreviated notation also for spaces of vector-valued or matrix-valued functions.

\section*{\S 1. Homogenization problem for a periodic elliptic operator
in $L_2(\mathbb{R}^d;\mathbb{\mathbb{C}}^n)$}

In this section, we describe the class of matrix elliptic
operators under consideration and formulate the results
on homogenization problem in $\mathbb{R}^d$ obtained in [BSu2,5].

\smallskip\noindent\textbf{1.1. Lattices in $\mathbb{R}^d$.}
Let $\mathbf{a}_1,\dots,\mathbf{a}_d \in \mathbb{R}^d$ be the basis in $\mathbb{R}^d$ that
generates the lattice $\Gamma$:
$$
\Gamma = \{ \mathbf{a} \in \mathbb{R}^d:\ \mathbf{a} = \sum_{j=1}^d \nu_j \mathbf{a}_j,
\ \ \nu_j \in \mathbb{Z}\},
$$
and let $\Omega$ be the elementary cell of $\Gamma$:
$$
\Omega := \{\mathbf{x} \in \mathbb{R}^d:\ \mathbf{x} =
\sum_{j=1}^d \tau_j \mathbf{a}_j,\ \ -\frac{1}{2} < \tau_j < \frac{1}{2} \}.
$$
We denote $|\Omega| = \text{meas}\, \Omega$.

The basis $\mathbf{b}_1,\dots, \mathbf{b}_d$ in $\mathbb{R}^d$ dual to $\mathbf{a}_1,\dots, \mathbf{a}_d$
is defined by the relations $\langle \mathbf{b}_i, \mathbf{a}_j\rangle  = 2 \pi \delta_{ij}$.
This basis generates the \textit{lattice} $\widetilde{\Gamma}$ \textit{dual to} $\Gamma$:
$$
\widetilde{\Gamma} = \{ \mathbf{b} \in \mathbb{R}^d:\ \mathbf{b} = \sum_{i=1}^d \rho_i
\mathbf{b}_i,\ \ \rho_i \in \mathbb{Z}\}.
$$
We introduce the \textit{central Brillouin zone}
$$
\widetilde{\Omega} = \{ \mathbf{k} \in \mathbb{R}^d:\  |\mathbf{k}| < |\mathbf{k}- \mathbf{b}|,
\ \ 0\ne \mathbf{b} \in \widetilde{\Gamma}\},
$$
which is a fundamental domain of $\widetilde{\Gamma}$.

Below, $\widetilde{H}^1(\Omega)$ stands for the subspace of all functions in
$H^1(\Omega)$ whose $\Gamma$-periodic extension to
$\mathbb{R}^d$ belongs to $H^1_{\text{loc}}(\mathbb{R}^d)$.
If $\varphi(\mathbf{x})$ is a $\Gamma$-periodic function in $\mathbb{R}^d$, we denote
$$
\varphi^\varepsilon(\mathbf{x}) := \varphi(\varepsilon^{-1}\mathbf{x}),\quad \varepsilon >0.
$$

\smallskip\noindent\textbf{1.2. The class of operators.}
In $L_2(\mathbb{R}^d;\mathbb{C}^n)$, we consider a second order DO
$\mathcal{A}_\varepsilon$ formally given by the differential expression
$$
\mathcal{A}_\varepsilon = b(\mathbf{D})^* g^\varepsilon(\mathbf{x}) b (\mathbf{D}),\ \ \varepsilon >0.
\eqno(1.1)
$$
Here $g(\mathbf{x})$ is a measurable $(m \times m)$-matrix-valued
function (in general, with complex entries). It is assumed that $g(\mathbf{x})$
is periodic with respect to the lattice $\Gamma$, bounded and uniformly positive definite.
Next, $b(\mathbf{D})$ is a homogeneous $(m\times n)$-matrix first order DO
with constant coefficients:
$$
b(\mathbf{D})=\sum_{l=1}^d b_l D_l.
\eqno(1.2)
$$
Here $b_l$ are constant matrices (in general, with complex entries).
The symbol $b(\boldsymbol{\xi}) = \sum_{l=1}^d b_l \xi_l$, $\boldsymbol{\xi} \in \mathbb{R}^d$,
is associated with the operator $b(\mathbf{D})$.
We assume that $m \ge n$ and that
$\textrm{rank} \, b (\boldsymbol{\xi}) = n, \ \forall \boldsymbol{\xi} \neq 0$.
This is equivalent to the inequalities
$$
\alpha_0 \mathbf{1}_n \leq b(\boldsymbol{\theta})^* b(\boldsymbol{\theta}) \leq \alpha_1 \mathbf{1}_n,
\quad \boldsymbol{\theta} \in \mathbb{S}^{d-1}, \quad 0 < \alpha_0 \leq \alpha_1 < \infty,
\eqno(1.3)
$$
with some positive constants $\alpha_0$ and $\alpha_1$.

The precise definition of the operator $\mathcal{A}_\varepsilon$ is given in terms
of the corresponding quadratic form
$$
a_\varepsilon[\mathbf{u},\mathbf{u}] = \int_{\mathbb{R}^d} \left\langle g^\varepsilon (\mathbf{x}) b (\mathbf{D}) \mathbf{u}, b(\mathbf{D}) \mathbf{u} \right\rangle \,
d \mathbf{x}, \quad \mathbf{u} \in H^1 (\mathbb{R}^{d}; \mathbb{C}^{n}).
$$
Under the above assumptions this form is closed in $L_2(\mathbb{R}^d;\mathbb{C}^n)$
and non-negative.
Using the Fourier transformation and (1.3), it is easy to check that
$$
c_0 \intop_{\mathbb{R}^d} |\mathbf{D} \mathbf{u}|^2\, d\mathbf{x}
\leq a_\varepsilon[\mathbf{u}, \mathbf{u}] \leq c_1 \intop_{\mathbb{R}^d}
|\mathbf{D} \mathbf{u}|^2\, d\mathbf{x},
\ \ \mathbf{u} \in H^1(\mathbb{R}^d;\mathbb{C}^n),
\eqno(1.4)
$$
where $c_0 = \alpha_0 \|g^{-1}\|^{-1}_{L_\infty}$, $c_1 = \alpha_1 \|g\|_{L_\infty}$.

The simplest example of the operator (1.1) is the scalar elliptic operator
$
\mathcal{A}_\varepsilon = -\text{div}\, g^\varepsilon(\mathbf{x}) \nabla
= \mathbf{D}^* g^\varepsilon(\mathbf{x})\mathbf{D}.
$
In this case we have $n=1$, $m=d$, $b(\mathbf{D})=\mathbf{D}$.
Obviously, (1.3) is true with $\alpha_0=\alpha_1=1$.
The operator of elasticity theory can be also written in the form (1.1)
with $n=d$, $m=d(d+1)/2$. These and other examples are considered in [BSu2] in detail.

\smallskip\noindent\textbf{1.3. The effective operator.}
In order to formulate the results, we need to introduce the effective operator $\mathcal{A}^0$.

Let an $(n\times m)$-matrix-valued function $\Lambda(\mathbf{x})$ be the (weak)
$\Gamma$-periodic solution of the problem
$$
b(\mathbf{D})^* g (\mathbf{x})\left( b(\mathbf{D}) \Lambda(\mathbf{x}) + \mathbf{1}_m \right) = 0,
\quad \int_{\Omega} \Lambda(\mathbf{x}) \, d \mathbf{x} = 0.
\eqno(1.5)
$$
In other words, for the columns $\mathbf{v}_j(\mathbf{x})$, $j=1,\dots,m,$
of the matrix $\Lambda(\mathbf{x})$ the following is true:
$\mathbf{v}_j \in \widetilde{H}^1(\Omega;\mathbb{C}^n)$, we have
$$
\int_\Omega \langle g(\mathbf{x}) (b(\mathbf{D})\mathbf{v}_j(\mathbf{x}) + \mathbf{e}_j),
b(\mathbf{D}) \boldsymbol{\eta}(\mathbf{x})\rangle\,d\mathbf{x} =0,
\quad \forall \boldsymbol{\eta} \in \widetilde{H}^1(\Omega;\mathbb{C}^n),
$$
and $\int_\Omega \mathbf{v}_j(\mathbf{x})\,d\mathbf{x}=0$.
Here $\mathbf{e}_1,\dots,\mathbf{e}_m$ is the standard orthonormal basis in $\mathbb{C}^m$.

The so-called \textit{effective matrix} $g^0$ of size $m\times m$ is defined as follows:
$$
g^0 = |\Omega|^{-1}
\int_{\Omega} g (\mathbf{x}) \left( b (\mathbf{D})
\Lambda (\mathbf{x}) + \mathbf{1}_m \right)\, d \mathbf{x}.
\eqno(1.6)
$$
It turns out that the matrix $g^0$ is positive definite.
The \textit{effective operator} $\mathcal{A}^0$ for the operator (1.1)
is given by the differential expression
$$
\mathcal{A}^0 = b (\mathbf{D})^* g^0 b (\mathbf{D})
$$
on the domain $H^2(\mathbb{R}^d;\mathbb{C}^n)$.

\smallskip\noindent\textbf{1.4. Properties of the effective matrix}.
The following properties of the effective matrix are proved in [BSu2, Chapter 3, Theorem 1.5].

\smallskip\noindent\textbf{Proposition 1.1.} \textit{The effective matrix} $g^0$
\textit{satisfies the following estimates}:
$$
\underline{g} \le g^0 \le \overline{g}.
\eqno(1.7)
$$
\textit{Here}
$$
\overline{g}= |\Omega|^{-1} \int_\Omega g(\mathbf{x})\,d\mathbf{x}, \quad
\underline{g}= \left(|\Omega|^{-1} \int_\Omega g(\mathbf{x})^{-1}\,d\mathbf{x}\right)^{-1}.
$$
\textit{If} $m = n$, \textit{then} $g^0$ \textit{coincides with} $\underline{g}$.

\smallskip
In homogenization theory, estimates (1.7) are well known for specific DO's as the Voight-Reuss bracketing.
Now we distinguish the cases where one of the inequalities in (1.7)
becomes an identity. The following statements were checked in [BSu2, Chapter 3, Propositions 1.6 and 1.7].

\smallskip\noindent\textbf{Proposition 1.2.} \textit{The identity} $g^0 = \overline{g}$
\textit{is equivalent to the relations}
$$
b(\mathbf{D})^* {\mathbf g}_k(\mathbf{x}) = 0, \ \ k = 1,\dots,m,
\eqno(1.8)
$$
\textit{where} ${\mathbf g}_k(\mathbf{x})$, $k = 1,\dots, m$, \textit{are the columns of
the matrix} $g(\mathbf{x})$.

\smallskip\noindent\textbf{Proposition 1.3.} \textit{The identity} $g^0 = \underline{g}$
\textit{is equivalent to the representations}
$$
 {\mathbf l}_k(\mathbf{x}) = {\mathbf l}_k^0 + b(\mathbf{D}) \mathbf{w}_k, \ \ {\mathbf l}_k^0 \in \mathbb{C}^m,
 \ \ \mathbf{w}_k \in \widetilde{H}^1(\Omega;\mathbb{C}^n),\ \ k = 1,\dots,m,
\eqno(1.9)
$$
 \textit{where} ${\mathbf l}_k(\mathbf{x})$, $k= 1,\dots, m$,
 \textit{are the columns of the matrix} $g(\mathbf{x})^{-1}$.

\smallskip
Obviously, (1.7) implies the following estimates for the norms of the matrices $g^0$ and $(g^0)^{-1}$:
$$
| g^0| \le \|g\|_{L_\infty},\ \ | (g^0)^{-1}| \le \|g^{-1}\|_{L_\infty}.
\eqno(1.10)
$$

\smallskip\noindent\textbf{1.5. The smoothing operator.}
We need an auxiliary smoothing operator $\Pi_\varepsilon$ acting in $L_2(\mathbb{R}^d;\mathbb{C}^m)$
and defined by the relation
$$
\left( \Pi_{\varepsilon} \mathbf{u} \right) (\mathbf{x}) =
(2 \pi)^{-d/2} \int_{\widetilde{\Omega}/\varepsilon}
e^{i \langle \mathbf{x}, \boldsymbol{\xi} \rangle} \widehat{\mathbf{u}} (\boldsymbol{\xi})
\, d \boldsymbol{\xi},
\eqno(1.11)
$$
where $\widehat{\mathbf{u}} (\boldsymbol{\xi})$ is the Fourier-image of $\mathbf{u}(\mathbf{x})$.
In other words, $\Pi_\varepsilon$ is the pseudodifferential operator with the symbol
$\chi_{\widetilde{\Omega}/\varepsilon} (\boldsymbol{\xi})$ which is the
indicator of the set $\widetilde{\Omega}/\varepsilon$.
Obviously, $\Pi_\varepsilon$ is the orthogonal projection in each space $H^s(\mathbb{R}^d;\mathbb{C}^m)$, $s \geq 0$.
Besides, $D^{\alpha} \Pi_{\varepsilon} \mathbf{u} = \Pi_{\varepsilon} D^{\alpha}\mathbf{u}$
for $\mathbf{u} \in H^s(\mathbb{R}^d; \mathbb{C}^m)$ and any multiindex $\alpha$ such that $|\alpha|\leq s$.

\smallskip\noindent\textbf{Proposition 1.4.}
\textit{For any} $\mathbf{u} \in H^1(\mathbb{R}^d;\mathbb{C}^m)$ \textit{we have}
$$
\| \Pi_\varepsilon \mathbf{u} - \mathbf{u} \|_{L_2(\mathbb{R}^d;\mathbb{C}^m)}
\leq \varepsilon r_0^{-1} \| \mathbf{D} \mathbf{u}\|_{L_2(\mathbb{R}^d)},
$$
\textit{where} $r_0$ \textit{is the radius of the ball inscribed in} $\text{clos}\,\widetilde{\Omega}$.

\smallskip\noindent\textbf{Proof.}
For $\boldsymbol{\xi} \in \mathbb{R}^d \setminus (\widetilde{\Omega}/\varepsilon)$ we have
$|\boldsymbol{\xi}| \geq r_0 \varepsilon^{-1}$. Hence,
$$
\begin{aligned}
&\| \Pi_\varepsilon \mathbf{u} - \mathbf{u} \|^2_{L_2(\mathbb{R}^d;\mathbb{C}^m)} =
\intop_{\mathbb{R}^d \setminus (\widetilde{\Omega}/\varepsilon)}
|\widehat{\mathbf{u}}(\boldsymbol{\xi})|^2\,d\boldsymbol{\xi}
\\
&\leq \varepsilon^2 r_0^{-2} \int_{\mathbb{R}^d}
|\boldsymbol{\xi}|^2 |\widehat{\mathbf{u}}(\boldsymbol{\xi})|^2\,d\boldsymbol{\xi}
=\varepsilon^2 r_0^{-2} \int_{\mathbb{R}^d} |\mathbf{D} \mathbf{u}(\mathbf{x})|^2 \,d\mathbf{x}. \ \ \bullet
\end{aligned}
$$

The following property was proved in [BSu5, Subsection 10.2].

\smallskip\noindent\textbf{Proposition 1.5.}
\textit{Let} $f(\mathbf{x})$ \textit{be a} $\Gamma$-\textit{periodic function in}
$\mathbb{R}^d$ \textit{such that}
$f \in L_2(\Omega)$. \textit{Let} $[f^\varepsilon]$ \textit{denote the operator of multiplication
by the function} $f(\varepsilon^{-1}\mathbf{x})$.
\textit{Then the operator} $[f^\varepsilon] \Pi_\varepsilon$ \textit{is continuous in}
$L_2(\mathbb{R}^d;\mathbb{C}^m)$, \textit{and}
$$
\| [f^\varepsilon] \Pi_\varepsilon\|_{L_2(\mathbb{R}^d;\mathbb{C}^m)\to L_2(\mathbb{R}^d;\mathbb{C}^m)}
\le |\Omega|^{-1/2} \|f\|_{L_2(\Omega)}.
$$

\smallskip\noindent\textbf{1.6. Results for homogenization problem in $\mathbb{R}^d$}.
Consider the following elliptic equation in $\mathbb{R}^d$:
$$
\mathcal{A}_\varepsilon \mathbf{u}_{\varepsilon} + \mathbf{u}_{\varepsilon} = \mathbf{F},
\eqno(1.12)
$$
where $\mathbf{F} \in L_2(\mathbb{R}^d;\mathbb{C}^n)$.
It is known that, as $\varepsilon \to 0$, the solution $\mathbf{u}_{\varepsilon}$ converges in
$L_2 (\mathbb{R}^d; \mathbb{C}^n)$ to the solution of the "homogenized" equation
$$
{\mathcal{A}}^0 \mathbf{u}_0 + \mathbf{u}_0 = \mathbf{F}.
\eqno(1.13)
$$

The following result was obtained in [BSu2, Chapter 4, Theorem 2.1].

\smallskip\noindent\textbf{Theorem 1.6.} \textit{Let} $\mathbf{u}_\varepsilon$
\textit{be the solution of the equation}
(1.12), \textit{and let} $\mathbf{u}_0$ \textit{be the solution of the equation} (1.13).
\textit{Then}
$$
 \| \mathbf{u}_{\varepsilon} - \mathbf{u}_0 \|_{L_2 (\mathbb{R}^{d}; \mathbb{C}^{n})}
 \leq C_1 \varepsilon  \| \mathbf{F} \|_{L_2 (\mathbb{R}^{d}; \mathbb{C}^{n})},
 \quad 0 < \varepsilon \leq 1,
 $$
\textit{or, in operator terms},
$$
\| (\mathcal{A}_{\varepsilon} + I)^{-1} - (\mathcal{A}^0 + I)^{-1}
\|_{L_2 (\mathbb{R}^{d}; \mathbb{C}^{n}) \rightarrow L_2 (\mathbb{R}^{d}; \mathbb{C}^{n})}
\leq \ C_1 \varepsilon, \quad 0 < \varepsilon \leq 1.
$$
\textit{The constant} $C_1$ \textit{depends only on the norms}
$\| g \|_{L_\infty}$, $\| g^{-1} \|_{L_\infty}$, \textit{the constants}
$\alpha_0$, $\alpha_1$ \textit{from} (1.3), \textit{and the parameters of the lattice} $\Gamma$.

\smallskip
In order to find approximation of the solution
$\mathbf{u}_\varepsilon$ in $H^1(\mathbb{R}^d;\mathbb{C}^n)$,
it is necessary to take the fist order corrector into account. We put
$$
K(\varepsilon) = [\Lambda^{\varepsilon}] \Pi_{\varepsilon} b (\mathbf{D}) (\mathcal{A}^0 +I)^{-1}.
\eqno(1.14)
$$
Here $[\Lambda^{\varepsilon}]$ is the operator of multiplication by the matrix-valued function
$\Lambda(\varepsilon^{-1}\mathbf{x})$, and
$\Pi_\varepsilon$ is the smoothing operator defined by (1.11).
The operator (1.14) is continuous from
$L_2(\mathbb{R}^d;\mathbb{C}^n)$ to $H^1(\mathbb{R}^d;\mathbb{C}^n)$.
This fact can be easily checked by using Proposition 1.5
and relation $\Lambda \in \widetilde{H}^1(\Omega)$.
Herewith, $\varepsilon \|K(\varepsilon)\|_{L_2 \to H^1} =O(1)$.

The "first order approximation" of the solution $\mathbf{u}_\varepsilon$ is given by
$$
\mathbf{v}_{\varepsilon} = \mathbf{u}_0 + \varepsilon \Lambda^{\varepsilon}
\Pi_{\varepsilon} b (\mathbf{D}) \mathbf{u}_0=
(\mathcal{A}^0+I)^{-1} \mathbf{F} + \varepsilon K(\varepsilon) \mathbf{F}.
\eqno(1.15)
$$
The following theorem was obtained in [BSu5, Theorem 10.6].

\smallskip\noindent\textbf{Theorem 1.7.} \textit{Let} $\mathbf{u}_\varepsilon$
\textit{be the solution of the equation} (1.12), \textit{and let} $\mathbf{u}_0$
\textit{be the solution of the equation} (1.13). \textit{Let}
$\mathbf{v}_\varepsilon$ \textit{be the function defined by} (1.15). \textit{Then}
$$
 \| \mathbf{u}_{\varepsilon} - \mathbf{v}_\varepsilon \|_{H^1 (\mathbb{R}^{d}; \mathbb{C}^{n})}
 \leq C_2 \varepsilon  \| \mathbf{F} \|_{L_2 (\mathbb{R}^{d}; \mathbb{C}^{n})},
 \quad 0 < \varepsilon \leq 1,
 \eqno(1.16)
 $$
\textit{or, in operator terms},
$$
\| (\mathcal{A}_{\varepsilon} + I)^{-1} - (\mathcal{A}^0 + I)^{-1} -
\varepsilon K (\varepsilon) \|_{L_2 (\mathbb{R}^{d}; \mathbb{C}^{n})
\to H^1 (\mathbb{R}^{d}; \mathbb{C}^{n})} \leq \ C_2 \varepsilon,
\quad 0 < \varepsilon \leq 1.
$$
\textit{The constant} $C_2$ \textit{depends only on}
$m, \alpha_0, \alpha_1, \| g \|_{L_\infty}, \| g^{-1} \|_{L_\infty}$,
\textit{and the parameters of the lattice} $\Gamma$.

\smallskip
Now we distinguish the case where the corrector is equal to zero.
Next statement follows from Theorem 1.7, Proposition 1.2 and equation (1.5).

\smallskip\noindent\textbf{Proposition 1.8.} \textit{If $g^0=\overline{g}$, i.~e.,
if relations} (1.8) \textit{are satisfied, then} $\Lambda=0$ \textit{and} $K(\varepsilon)=0$.
\textit{Then we have}
$$
 \| \mathbf{u}_{\varepsilon} - \mathbf{u}_0 \|_{H^1 (\mathbb{R}^{d}; \mathbb{C}^{n})}
 \leq C_2 \varepsilon  \| \mathbf{F} \|_{L_2 (\mathbb{R}^{d}; \mathbb{C}^{n})},
 \quad 0 < \varepsilon \leq 1.
 $$

\smallskip

It turns out that under some assumptions on the solution of the problem (1.5)
the smoothing operator $\Pi_{\varepsilon}$ in the corrector (1.14) can be removed
(replaced by the identity).

\smallskip\noindent\textbf{Condition 1.9.}
\textit{Suppose that the} $\Gamma$-\textit{periodic solution}
$\Lambda(\mathbf{x})$ \textit{of the problem} (1.5) \textit{is bounded}:
$
\Lambda \in L_\infty.
$

\smallskip
We put
$$
K^0 (\varepsilon) = [\Lambda^{\varepsilon}] b (\mathbf{D}) (\mathcal{A}^0 + I)^{-1}.
$$
In [BSu5], it was shown that under Condition 1.9
the operator $K^0 (\varepsilon)$ is a continuous mapping of
$L_2 (\mathbb{R}^{d}; \mathbb{C}^{n})$ into $H^1 (\mathbb{R}^{d}; \mathbb{C}^{n})$.
(It is also easy to deduce this fact from Corollary 2.4 proved below.)

Instead of (1.15), we consider another approximation of the solution $\mathbf{u}_{\varepsilon}$:
$$
\check{\mathbf{v}}_{\varepsilon} = \mathbf{u}_0 + \varepsilon \Lambda^{\varepsilon} b (\mathbf{D}) \mathbf{u}_0
= (\mathcal{A}^0 + I)^{-1} \mathbf{F} + \varepsilon K^0 (\varepsilon) \mathbf{F}.
\eqno(1.17)
$$

The following result was obtained in [BSu5, Theorem 10.8].

\smallskip\noindent\textbf{Theorem 1.10.} \textit{Suppose that Condition} 1.9
\textit{is satisfied. Let} $\mathbf{u}_\varepsilon$ \textit{be the solution of the equation}
(1.12), \textit{and let} $\mathbf{u}_0$ \textit{be the solution of the equation} (1.13).
\textit{Let} $\check{\mathbf{v}}_\varepsilon$ \textit{be the function defined by} (1.17).
\textit{Then we have}
$$
 \| \mathbf{u}_{\varepsilon} - \check{\mathbf{v}}_\varepsilon \|_{H^1 (\mathbb{R}^{d}; \mathbb{C}^{n})}
 \leq C_3 \varepsilon  \| \mathbf{F} \|_{L_2 (\mathbb{R}^{d}; \mathbb{C}^{n})},
 \quad 0 < \varepsilon \leq 1,
 $$
\textit{or, in operator terms,}
$$
\| (\mathcal{A}_{\varepsilon} + I)^{-1} - (\mathcal{A}^0 + I)^{-1} -
\varepsilon K^0 (\varepsilon) \|_{L_2 (\mathbb{R}^{d}; \mathbb{C}^{n})
\to H^1 (\mathbb{R}^{d}; \mathbb{C}^{n})} \leq \ C_3 \varepsilon,
\quad 0 < \varepsilon \leq 1.
$$
\textit{The constant} $C_3$ \textit{depends only on}
$m, d, \alpha_0, \alpha_1, \| g \|_{L_\infty}, \| g^{-1} \|_{L_\infty}$,
\textit{the parameters of the lattice} $\Gamma$,
\textit{and the norm} $\|\Lambda\|_{L_\infty}$.

\smallskip
In some cases Condition 1.9 is valid automatically.
The following statement was checked in [BSu5, Lemma 8.7].

\smallskip\noindent\textbf{Proposition 1.11.} \textit{Condition} 1.9 \textit{is a fortiori valid if
at least one of the following assumptions is satisfied}:

$1^\circ$. \textit{dimension does not exceed two, i.~e.} $d \leq 2$;

$2^\circ$. \textit{the operator acts in} $L_2(\mathbb{R}^d)$, $d\ge 1$, \textit{and has the form}
$\mathcal{A}_\varepsilon = \mathbf{D}^* g^\varepsilon(\mathbf{x}) \mathbf{D}$,
\textit{where} $g(\mathbf{x})$ \textit{is a matrix with real entries};

$3^\circ$. \textit{dimension is arbitrary and} $g^0 = \underline{g}$, \textit{i.~e., relations}
(1.9) \textit{are satisfied}.

\smallskip
Note that Condition 1.9 can be also ensured by the assumption that
the matrix $g(\mathbf{x})$ is sufficiently smooth.

\section*{\S 2. Properties of the matrix-valued function $\Lambda$}

The following statement is proved by analogy with the proof of
Lemma 8.3 from [BSu5].

\smallskip\noindent\textbf{Lemma 2.1.} \textit{Let}
$\Lambda (\mathbf{x})$ \textit{be the} $\Gamma$-\textit{periodic solution of the problem} (1.5).
\textit{Then for any function} $u \in C_0^{\infty} (\mathbb{R}^d)$ \textit{we have}
$$
\int_{\mathbb{R}^d} |\mathbf{D} \Lambda(\mathbf{x})|^2 |u|^2 \, d \mathbf{x}
\leq \beta_1 \| u \|^2_{L_2 (\mathbb{R}^{d})} +
\beta_2 \int_{\mathbb{R}^d} |\Lambda(\mathbf{x})|^2 |\mathbf{D} u|^2 \, d \mathbf{x}.
\eqno(2.1)
$$
\textit{The constants} $\beta_1$ \textit{and} $\beta_2$
\textit{are defined below in} (2.12) \textit{and depend only on} $m$, $d$, $\alpha_0$,
$\alpha_1$, $\| g \|_{L_\infty}$, \textit{and} $\| g^{-1} \|_{L_\infty}$.

\smallskip\noindent\textbf{Proof.}
Let $\mathbf{v}_j(\mathbf{x})$, $j=1,\dots,m$, be the columns of the matrix $\Lambda (\mathbf{x})$.
By (1.5), for any function $\boldsymbol{\eta} \in H^1 (\mathbb{R}^d; \mathbb{C}^{n})$ such that
$\boldsymbol{\eta} (\mathbf{x}) = 0$ for $|\mathbf{x}| > r$ (with some $r > 0$)
we have
$$
\int_{\mathbb{R}^d} \left\langle g (\mathbf{x}) \left( b (\mathbf{D})
\mathbf{v}_j (\mathbf{x}) + \mathbf{e}_j \right), b(\mathbf{D}) \boldsymbol{\eta} (\mathbf{x})
\right\rangle \, d \mathbf{x} = 0.
\eqno(2.2)
$$

Let $u \in C_0^\infty(\mathbb{R}^d)$. We put
$\boldsymbol{\eta} (\mathbf{x}) = \mathbf{v}_j (\mathbf{x}) |u (\mathbf{x})|^2$.
By (1.2),
$$
b(\mathbf{D}) \boldsymbol{\eta} (\mathbf{x}) =
\left( b(\mathbf{D}) \mathbf{v}_j (\mathbf{x}) \right) | u (\mathbf{x}) |^2 +
\sum_{l=1}^d b_l \mathbf{v}_j (\mathbf{x}) D_l | u (\mathbf{x}) |^2.
\eqno(2.3)
$$
Substituting (2.3) in (2.2), we obtain
$$
\begin{aligned}
&\int_{\mathbb{R}^d} \left\langle  g (\mathbf{x})
\left( b (\mathbf{D}) \mathbf{v}_j (\mathbf{x}) + \mathbf{e}_j \right),
b(\mathbf{D}) \mathbf{v}_j (\mathbf{x}) \right\rangle |u|^2 \, d \mathbf{x}
\\
&+ \int_{\mathbb{R}^d} \sum_{l=1}^d \langle g (\mathbf{x}) \left( b(\mathbf{D})\mathbf{v}_j (\mathbf{x}) + \mathbf{e}_j
\right),b_l \mathbf{v}_j (\mathbf{x})\rangle
\left( D_l u \, \overline{u} + u \, D_l \overline{u} \right) \, d \mathbf{x} = 0.
\end{aligned}
$$
Hence,
$$
\begin{aligned}
J:&=\int_{\mathbb{R}^d} \left| g^{1/2} b(\mathbf{D}) \mathbf{v}_j \right|^2 |u|^2 \, d \mathbf{x}
= - \int_{\mathbb{R}^d} \left\langle g^{1/2} \mathbf{e}_j, g^{1/2} b(\mathbf{D})\mathbf{v}_j
\right\rangle |u|^2 \, d \mathbf{x}
\\
&- \int_{\mathbb{R}^d} \sum_{l=1}^d \langle g^{1/2} b(\mathbf{D}) \mathbf{v}_j,
g^{1/2} b_l \mathbf{v}_j \rangle \left( D_l u \, \overline{u} + u \, D_l \overline{u} \right)
\, d \mathbf{x}
\\
&- \int_{\mathbb{R}^d} \sum_{l=1}^d \langle g \mathbf{e}_j, b_l \mathbf{v}_j \rangle
\left( D_l u \, \overline{u} + u \, D_l \overline{u} \right) \, d \mathbf{x}.
\end{aligned}
\eqno(2.4)
$$
Denote the summands on the right by $J_1$, $J_2$, $J_3$.
The first term $J_1$ can be estimated as follows:
$$
| J_1 | \leq \int_{\mathbb{R}^d} \left( \left| g^{1/2} \mathbf{e}_j \right|^2 +
\frac{1}{4} \left| g^{1/2} b(\mathbf{D}) \mathbf{v}_j \right|^2 \right) |u|^2 \, d \mathbf{x}
\leq  \| g \|_{L_\infty} \|u\|^2_{L_2(\mathbb{R}^d)} + \frac{1}{4} J.
\eqno(2.5)
$$

Next, from (1.3) it follows that
$$
|b_l| \leq \alpha_1^{1/2}, \quad l = 1, \ldots, d.
\eqno(2.6)
$$
Taking (2.6) into account, we estimate the second term $J_2$:
$$
\begin{aligned}
| J_2 | &\leq 2 \int_{\mathbb{R}^d} \left| g^{1/2} \ b(\mathbf{D}) \mathbf{v}_j \right| \,|u|
\left(\sum_{l=1}^d \left| g^{1/2} b_l \mathbf{v}_j \right| |D_l u| \right) \, d \mathbf{x}
\\
&\leq \frac{1}{4} J  + 4 d \alpha_1 \| g \|_{L_\infty}
\int_{\mathbb{R}^d} | \mathbf{v}_j |^2 |\mathbf{D} u|^2 \, d \mathbf{x}.
\end{aligned}
\eqno(2.7)
$$
Finally, the term $J_3$ satisfies the estimate
$$
\begin{aligned}
| J_3 | &\leq 2 \int_{\mathbb{R}^d} \left| g \mathbf{e}_j \right| \,|u|
\left( \sum_{l=1}^d \left| b_l \mathbf{v}_j \right| |D_l u| \right) d \mathbf{x}
\\
&\leq \| g \|_{L_\infty} \|u\|^2 _{L_2(\mathbb{R}^d)} +
d \alpha_1 \| g \|_{L_\infty}
\int_{\mathbb{R}^d} | \mathbf{v}_j |^2 |\mathbf{D} u|^2 \, d \mathbf{x}.
\end{aligned}
\eqno(2.8)
$$

Combining (2.4), (2.5), (2.7), and (2.8), we obtain
$$
\frac{1}{2} J \leq 2 \| g \|_{L_\infty} \|u\|^2_{L_2(\mathbb{R}^d)} +
5 d \alpha_1 \| g \|_{L_\infty}
\int_{\mathbb{R}^d} | \mathbf{v}_j |^2 |\mathbf{D} u|^2 \, d \mathbf{x}.
\eqno(2.9)
$$

Now, we show how the required estimate can be deduced from (2.9).
By the Fourier transformation, it follows from the lower inequality (1.3) that
$$
 \int_{\mathbb{R}^d} \left| \mathbf{D} (\mathbf{v}_j u) \right|^2 \, d \mathbf{x} \leq
 \alpha_0^{-1} \int_{\mathbb{R}^d} \left| b(\mathbf{D}) (\mathbf{v}_j u) \right|^2 \ d \mathbf{x}.
$$
By (1.2),
$$
b(\mathbf{D}) (\mathbf{v}_j u) = (b(\mathbf{D})\mathbf{v}_j) u + \sum_{l=1}^d b_l \mathbf{v}_j D_l u.
$$
Then, taking (2.6) and the expression for $J$ (see (2.4)) into account, we have
$$
\begin{aligned}
\int_{\mathbb{R}^d} \left| \mathbf{D} (\mathbf{v}_j u) \right|^2 \, d \mathbf{x} &\leq
2 \alpha_0^{-1} \int_{\mathbb{R}^d} \left| b(\mathbf{D}) \mathbf{v}_j \right|^2 |u|^2 \,
d \mathbf{x} + 2 \alpha_0^{-1} \alpha_1 d \int_{\mathbb{R}^d} |\mathbf{v}_j|^2 |\mathbf{D} u|^2 \, d \mathbf{x}
\\
&\leq 2 \alpha_0^{-1} \| g^{-1} \|_{L_\infty} J +
2 \alpha_0^{-1} \alpha_1 d \int_{\mathbb{R}^d} |\mathbf{v}_j|^2 |\mathbf{D} u|^2 \, d \mathbf{x}.
\end{aligned}
\eqno(2.10)
$$
Obviously,
$$
\int_{\mathbb{R}^d} \left| \mathbf{D} \mathbf{v}_j \right|^2 |u|^2 \, d \mathbf{x}
\leq 2 \int_{\mathbb{R}^d} \left| \mathbf{D} (\mathbf{v}_j u) \right|^2 \, d \mathbf{x}
+ 2 \int_{\mathbb{R}^d} | \mathbf{v}_j |^2 |\mathbf{D} u|^2 \, d \mathbf{x}.
\eqno(2.11)
$$
Relations (2.9)--(2.11) imply that
$$
\begin{aligned}
&\int_{\mathbb{R}^d} \left| \mathbf{D} \mathbf{v}_j \right|^2 |u|^2 \, d \mathbf{x}
\leq  16 \alpha_0^{-1} \| g^{-1} \|_{L_\infty} \| g \|_{L_\infty} \| u \|^2_{L_2(\mathbb{R}^d)}
\\
&+ 2 \left(1 + 2 d \alpha_0^{-1} \alpha_1 + 20 d \alpha_0^{-1} \alpha_1
\| g^{-1} \|_{L_\infty} \| g \|_{L_\infty} \right)
\int_{\mathbb{R}^d} | \mathbf{v}_j |^2 |\mathbf{D} u|^2 \, d \mathbf{x}.
\end{aligned}
$$
Summing up over $j$, we arrive at estimate (2.1) with
$$
\begin{aligned}
\beta_1 &= 16 m \alpha_0^{-1} \| g^{-1} \|_{L_\infty} \| g \|_{L_\infty},
\\
\beta_2 &=2 \left(1 + 2 d \alpha_0^{-1} \alpha_1 + 20 d \alpha_0^{-1} \alpha_1
\| g^{-1} \|_{L_\infty} \| g \|_{L_\infty} \right).
\ \ \bullet
\end{aligned}
\eqno(2.12)
$$

\smallskip\noindent\textbf{Corollary 2.2.} \textit{Under Condition} 1.9
\textit{for any} $u \in H^1 (\mathbb{R}^d)$ \textit{we have}
$$
\int_{\mathbb{R}^d} |\mathbf{D} \Lambda(\mathbf{x})|^2 |u|^2 \, d \mathbf{x}
\leq \beta_1 \| u \|^2_{L_2 (\mathbb{R}^{d})} +
\beta_2 \|\Lambda\|^2_{L_\infty} \int_{\mathbb{R}^d} |\mathbf{D} u|^2 \, d \mathbf{x}.
\eqno(2.13)
$$

\smallskip\noindent\textbf{Proof.}
Indeed, the second integral in the right-hand side of (2.1)
can be estimated by $\| \Lambda \|^2_{L_\infty} \int_{\mathbb{R}^d} | \mathbf{D} u |^2 \, d \mathbf{x}$.
Then (2.13) is valid for any $u\in C_0^\infty(\mathbb{R}^d)$. By continuity
inequality (2.13) is extended from the dense set
$C_0^{\infty} (\mathbb{R}^d)$ to the whole $H^1 (\mathbb{R}^d)$. $\bullet$

Next statement follows from Lemma 2.1 by the scaling transformation.

\smallskip\noindent\textbf{Lemma 2.3.} \textit{Under the assumptions of Lemma}
2.1 \textit{we have}
$$
\int_{\mathbb{R}^d} \left| \left( \mathbf{D} \Lambda \right)^{\varepsilon} (\mathbf{x}) \right|^2
|u(\mathbf{x})|^2 \, d \mathbf{x}  \leq \beta_1  \| u \|^2_{L_2 (\mathbb{R}^{d})} +
\beta_2  \varepsilon^2 \int_{\mathbb{R}^d}
\left| \Lambda^{\varepsilon} (\mathbf{x}) \right|^2 |\mathbf{D} u|^2 \, d \mathbf{x}.
$$

\smallskip\noindent\textbf{Proof.}
By the changes $\mathbf{y} = \varepsilon^{-1} \mathbf{x}$
and $u(\mathbf{x}) = v (\mathbf{y})$, from (2.1) it follows that
$$
\begin{aligned}
&\int_{\mathbb{R}^d} \left| \left( \mathbf{D} \Lambda \right)
( \varepsilon^{-1}\mathbf{x}) \right|^2
|u (\mathbf{x})|^2 \, d \mathbf{x} = \int_{\mathbb{R}^d} \left| \left( \mathbf{D} \Lambda \right)
\left( \mathbf{y} \right) \right|^2 |v (\mathbf{y})|^2 \, \varepsilon^d \, d \mathbf{y}
\\
&\leq \beta_1 \int_{\mathbb{R}^d} |v (\mathbf{y})|^2 \varepsilon^d d \mathbf{y}
+ \beta_2 \int_{\mathbb{R}^d} \left| \Lambda (\mathbf{y}) \right|^2
|\mathbf{D}_y v (\mathbf{y})|^2 \varepsilon^d \, d \mathbf{y}
\\
&=\beta_1 \int_{\mathbb{R}^d} |u (\mathbf{x})|^2 \, d \mathbf{x} +
\beta_2 \varepsilon^2 \int_{\mathbb{R}^d}
\left| \Lambda \left( \varepsilon^{-1} \mathbf{x}\right) \right|^2
|\mathbf{D}_\mathbf{x} u (\mathbf{x})|^2 \, d \mathbf{x}. \quad \bullet
\end{aligned}
$$

\smallskip\noindent\textbf{Corollary 2.4.} \textit{Under Condition} 1.9
\textit{for any} $u \in H^1 (\mathbb{R}^d)$ \textit{we have}
$$
\int_{\mathbb{R}^d} |(\mathbf{D} \Lambda)^\varepsilon(\mathbf{x})|^2 |u|^2 \, d \mathbf{x}
\leq \beta_1 \| u \|^2_{L_2 (\mathbb{R}^{d})} +
\beta_2 \|\Lambda\|^2_{L_\infty} \varepsilon^2 \int_{\mathbb{R}^d} |\mathbf{D} u|^2 \, d \mathbf{x}.
$$

\smallskip
In conclusion of this section, we give two estimates for the matrix-valued function $\Lambda$ obtained in
[BSu4, (6.28) and Subsection 7.3]:
$$
\| \Lambda \|_{L_2(\Omega)} \leq |\Omega|^{1/2} m^{1/2} (2r_0)^{-1} \alpha_0^{-1/2}
\|g\|_{L_\infty}^{1/2} \|g^{-1} \|_{L_\infty}^{1/2},
\eqno(2.14)
$$
$$
\| \mathbf{D} \Lambda \|_{L_2(\Omega)} \leq |\Omega|^{1/2} m^{1/2} \alpha_0^{-1/2}
\|g\|_{L_\infty}^{1/2} \|g^{-1} \|_{L_\infty}^{1/2}.
\eqno(2.15)
$$

\section*{\S 3. Smoothing in Steklov's sense. One more result for homogenization problem in $\mathbb{R}^d$}

In [Zh2, ZhPas], smoothing in Steklov's sense was used
instead of the smoothing operator (1.11).
It turns out that smoothing in Steklov's sense is more convenient for the study of homogenization problem
in a bounded domain. In this section, we show that for the problem in
$\mathbb{R}^d$ both variants are possible, i.~e., Theorem 1.7 remains true if in the corrector (1.14)
the operator $\Pi_\varepsilon$ is replaced by the operator smoothing in Steklov's sense.

\smallskip\noindent\textbf{3.1. Smoothing in Steklov's sense.}
In $L_2(\mathbb{R}^d;\mathbb{C}^m)$, we consider the operator $S_\varepsilon$ defined by
$$
(S_\varepsilon \mathbf{u})(\mathbf{x}) = |\Omega|^{-1}
\int_\Omega \mathbf{u}(\mathbf{x} - \varepsilon \mathbf{z})\, d\mathbf{z}
\eqno(3.1)
$$
and called the \textit{operator smoothing in Steklov's sense}.
It is easy to check that $\|S_\varepsilon\|_{L_2(\mathbb{R}^d) \to L_2(\mathbb{R}^d)} \leq 1$.
Obviously, $D^\alpha S_\varepsilon \mathbf{u} = S_\varepsilon D^\alpha \mathbf{u}$
for $\mathbf{u} \in H^s(\mathbb{R}^d;\mathbb{C}^m)$ and any multiindex $\alpha$ such that $|\alpha| \leq s$.

We need some properties of the operator (3.1), cf. [ZhPas, Lemmas 1.1 and 1.2].

\smallskip\noindent\textbf{Proposition 3.1.}
\textit{For any} $\mathbf{u} \in H^1(\mathbb{R}^d;\mathbb{C}^m)$ \textit{we have}
$$
\| S_\varepsilon \mathbf{u} - \mathbf{u}\|_{L_2(\mathbb{R}^d;\mathbb{C}^m)} \leq \varepsilon r_1 \| \mathbf{D} \mathbf{u} \|_{L_2(\mathbb{R}^d)},
\eqno(3.2)
$$
\textit{where} $2 r_1 =\text{diam}\,\Omega$.

\smallskip\noindent\textbf{Proof.} By the Cauchy inequality,
$$
\begin{aligned}
\| S_\varepsilon \mathbf{u} - \mathbf{u}\|_{L_2(\mathbb{R}^d;\mathbb{C}^m)}^2
&= \int_{\mathbb{R}^d} d\mathbf{x} \left| |\Omega|^{-1} \int_\Omega (\mathbf{u}(\mathbf{x} - \varepsilon \mathbf{z}) - \mathbf{u}(\mathbf{x}))\,d\mathbf{z} \right|^2
\\
&\leq |\Omega|^{-1} \int_{\mathbb{R}^d}d\mathbf{x} \int_\Omega|\mathbf{u}(\mathbf{x} - \varepsilon \mathbf{z}) - \mathbf{u}(\mathbf{x})|^2\,d\mathbf{z}.
\end{aligned}
\eqno(3.3)
$$
Using the Fourier transformation, we obtain
$$
\begin{aligned}
\int_{\mathbb{R}^d} |\mathbf{u}(\mathbf{x} - \varepsilon \mathbf{z}) - \mathbf{u}(\mathbf{x})|^2\,d\mathbf{x}=
\int_{\mathbb{R}^d} \left| \exp(-i \varepsilon \langle \mathbf{z},\boldsymbol{\xi}\rangle) -1 \right|^2
|\widehat{\mathbf{u}}(\boldsymbol{\xi})|^2 \,d\boldsymbol{\xi}
\\
\leq \varepsilon^2 |\mathbf{z}|^2 \int_{\mathbb{R}^d} |\boldsymbol{\xi}|^2|\widehat{\mathbf{u}}(\boldsymbol{\xi})|^2 \,d\boldsymbol{\xi}
= \varepsilon^2 |\mathbf{z}|^2 \int_{\mathbb{R}^d} |\mathbf{D} \mathbf{u}(\mathbf{x})|^2 \,d\mathbf{x}.
\end{aligned}
$$
Integrating this inequality over $\mathbf{z} \in \Omega$, we conclude that
$$
\int_\Omega d \mathbf{z} \int_{\mathbb{R}^d} |\mathbf{u}(\mathbf{x} - \varepsilon \mathbf{z}) - \mathbf{u}(\mathbf{x})|^2\,d\mathbf{x}
\leq \varepsilon^2 r_1^2 |\Omega| \int_{\mathbb{R}^d} | \mathbf{D} \mathbf{u} (\mathbf{x})|^2\,d\mathbf{x}.
$$
Together with (3.3) this implies (3.2). $\bullet$

\smallskip\noindent\textbf{Proposition 3.2.} \textit{Let} $f(\mathbf{x})$
\textit{be a} $\Gamma$-\textit{periodic function in} $\mathbb{R}^d$ \textit{such that} $f \in L_2(\Omega)$.
\textit{Then the operator} $[f^\varepsilon]S_\varepsilon$
\textit{is continuous in} $L_2(\mathbb{R}^d;\mathbb{C}^m)$, \textit{and}
$$
\| [f^\varepsilon]S_\varepsilon \|_{L_2(\mathbb{R}^d;\mathbb{C}^m)\to L_2(\mathbb{R}^d;\mathbb{C}^m)}
\leq |\Omega|^{-1/2} \| f \|_{L_2(\Omega)}.
$$

\smallskip\noindent\textbf{Proof.}
By the Cauchy inequality and the change of variables, from (3.1) it follows that
$$
\begin{aligned}
&\int_{\mathbb{R}^d} |f^\varepsilon(\mathbf{x}) (S_\varepsilon \mathbf{u})(\mathbf{x})|^2 \,d\mathbf{x} \leq
|\Omega|^{-1} \int_{\mathbb{R}^d} d\mathbf{x}\,|f(\varepsilon^{-1}\mathbf{x})|^2 \int_\Omega |\mathbf{u}(\mathbf{x}-\varepsilon \mathbf{z})|^2 \,d\mathbf{z}
\\
&=
|\Omega|^{-1} \int_{\mathbb{R}^d} d\mathbf{y}\int_\Omega |f(\varepsilon^{-1}\mathbf{y} + \mathbf{z})|^2 |\mathbf{u}(\mathbf{y})|^2 \,d\mathbf{z} =
|\Omega|^{-1} \|f\|^2_{L_2(\Omega)} \|\mathbf{u}\|^2_{L_2(\mathbb{R}^d)}.\ \bullet
\end{aligned}
$$

\smallskip\noindent\textbf{3.2.} We put
$$
\widetilde{K}(\varepsilon) = [\Lambda^\varepsilon] S_\varepsilon b(\mathbf{D}) (\mathcal{A}^0 +I)^{-1}.
\eqno(3.4)
$$
The operator (3.4) is continuous from $L_2(\mathbb{R}^d;\mathbb{C}^n)$ to $H^1(\mathbb{R}^d;\mathbb{C}^n)$.
Indeed, the operator $b(\mathbf{D}) (\mathcal{A}^0 +I)^{-1}$
is a continuous mapping of $L_2(\mathbb{R}^d;\mathbb{C}^n)$ into $H^1(\mathbb{R}^d;\mathbb{C}^m)$.
Using Proposition 3.2 and relation $\Lambda \in \widetilde{H}^1(\Omega)$, it is easy to
check that the operator $[\Lambda^\varepsilon] S_\varepsilon$ is continuous from
$H^1(\mathbb{R}^d;\mathbb{C}^m)$ to $H^1(\mathbb{R}^d;\mathbb{C}^n)$.

Let $\mathbf{u}_\varepsilon$ be the solution of the equation (1.12).
Instead of (1.15) we consider another fist order approximation of $\mathbf{u}_\varepsilon$:
$$
\widetilde{\mathbf{v}}_\varepsilon = \mathbf{u}_0 +
\varepsilon \Lambda^\varepsilon S_\varepsilon b(\mathbf{D}) \mathbf{u}_0 =
(\mathcal{A}^0 +I)^{-1} \mathbf{F} + \varepsilon \widetilde{K}(\varepsilon) \mathbf{F}.
\eqno(3.5)
$$
Along with Theorem 1.7, the following result takes place.

\smallskip\noindent\textbf{Theorem 3.3.} \textit{Let} $\mathbf{u}_\varepsilon$ \textit{be the solution
of the equation} (1.12), \textit{and let} $\mathbf{u}_0$ \textit{be the solution of the equation}
(1.13). \textit{Let} $\widetilde{\mathbf{v}}_\varepsilon$ \textit{be the function defined by} (3.5).
\textit{Then}
$$
\| \mathbf{u}_\varepsilon - \widetilde{\mathbf{v}}_\varepsilon \|_{H^1(\mathbb{R}^d;\mathbb{C}^n)} \leq \widetilde{C}_2 \varepsilon \| \mathbf{F} \|_{L_2(\mathbb{R}^d;\mathbb{C}^n)},
\ \ 0< \varepsilon \leq 1,
\eqno(3.6)
$$
\textit{or, in operator terms},
$$
\|(\mathcal{A}_\varepsilon +I)^{-1} - (\mathcal{A}^0 +I)^{-1} -
\varepsilon \widetilde{K}(\varepsilon) \|_{L_2(\mathbb{R}^d;\mathbb{C}^n) \to H^1(\mathbb{R}^d;\mathbb{C}^n)}
\leq \widetilde{C}_2 \varepsilon,\ \ 0< \varepsilon \leq 1.
$$
\textit{The constant} $\widetilde{C}_2$ \textit{depends only on} $m$, $d$, $\alpha_0$, $\alpha_1$,
$\|g\|_{L_\infty}$, $\|g^{-1}\|_{L_\infty}$, \textit{and the parameters of the lattice} $\Gamma$.

\smallskip
Theorem 3.3 will be deduced from Theorem 1.7.

\smallskip\noindent\textbf{Lemma 3.4.} \textit{For any} $\mathbf{u} \in H^2(\mathbb{R}^d;\mathbb{C}^n)$
\textit{we have}
$$
\begin{aligned}
\int_{\mathbb{R}^d}& |(\mathbf{D} \Lambda)^\varepsilon|^2 |(\Pi_\varepsilon - S_\varepsilon) b(\mathbf{D})\mathbf{u}|^2 \,d\mathbf{x} \leq
\beta_1 \int_{\mathbb{R}^d}|(\Pi_\varepsilon - S_\varepsilon) b(\mathbf{D})\mathbf{u}|^2 \,d\mathbf{x}
\\
&+ \beta_2 \varepsilon^2 \sum_{j=1}^d
\int_{\mathbb{R}^d} |\Lambda^\varepsilon|^2 |(\Pi_\varepsilon - S_\varepsilon) b(\mathbf{D}) \partial_j \mathbf{u}|^2\,d\mathbf{x}.
\end{aligned}
\eqno(3.7)
$$

\smallskip\noindent\textbf{Proof.}
By Propositions 1.5 and 3.2 and relation $\Lambda \in \widetilde{H}^1(\Omega)$,
all the terms in (3.7) are continuous functionals of $\mathbf{u}$
in the norm of $H^2(\mathbb{R}^d;\mathbb{C}^n)$.
Since $C_0^\infty(\mathbb{R}^d;\mathbb{C}^n)$ is dense in $H^2(\mathbb{R}^d;\mathbb{C}^n)$,
it suffices to check (3.7) for $\mathbf{u} \in C_0^\infty(\mathbb{R}^d;\mathbb{C}^n)$.

We fix a function $\zeta \in C^\infty(\mathbb{R}_+)$ such that $0\leq \zeta(t)\leq 1$,
$\zeta(t)=1$ for $0\leq t \leq 1$, and $\zeta(t)=0$ for $t\geq 2$.
We put $\zeta_R(\mathbf{x})= \zeta(R^{-1}|\mathbf{x}|)$, $\mathbf{x} \in \mathbb{R}^d$,
$R>0$. Let $\mathbf{u} \in C_0^\infty(\mathbb{R}^d;\mathbb{C}^n)$.
Then $\zeta_R (\Pi_\varepsilon - S_\varepsilon) b(\mathbf{D})\mathbf{u}
\in C_0^\infty(\mathbb{R}^d;\mathbb{C}^m)$ and, by Lemma 2.3, we have
$$
\begin{aligned}
\int_{\mathbb{R}^d} |(\mathbf{D} \Lambda)^\varepsilon|^2
|\zeta_R(\Pi_\varepsilon - S_\varepsilon) b(\mathbf{D})\mathbf{u}|^2 \,d\mathbf{x} \leq
\beta_1 \int_{\mathbb{R}^d}|\zeta_R(\Pi_\varepsilon - S_\varepsilon)
b(\mathbf{D})\mathbf{u}|^2 \,d\mathbf{x}
\\
+ \beta_2 \varepsilon^2 \sum_{j=1}^d
\int_{\mathbb{R}^d} |\Lambda^\varepsilon|^2
|(\partial_j \zeta_R)(\Pi_\varepsilon - S_\varepsilon) b(\mathbf{D})\mathbf{u} +
\zeta_R(\Pi_\varepsilon - S_\varepsilon) b(\mathbf{D}) \partial_j \mathbf{u}|^2\,d\mathbf{x}.
\end{aligned}
$$
Take into account that $\max\,|\partial_j \zeta_R| \leq c R^{-1}$. Then
(3.7) follows from the last inequality by the limit procedure as $R\to \infty$,
by the Lebesgue Theorem.  $\ \bullet$

From Proposition 1.5 and estimate (2.14) it follows that
$$
\begin{aligned}
\|[\Lambda^\varepsilon] \Pi_\varepsilon \|_{L_2(\mathbb{R}^d;\mathbb{C}^m) \to L_2(\mathbb{R}^d;\mathbb{C}^n)}
\leq |\Omega|^{-1/2} \| \Lambda \|_{L_2(\Omega)}
\\
\leq m^{1/2} (2r_0)^{-1} \alpha_0^{-1/2} \|g\|_{L_\infty}^{1/2} \|g^{-1}\|_{L_\infty}^{1/2}=: M.
\end{aligned}
\eqno(3.8)
$$
Similarly, Proposition 3.2 implies that
$$
\|[\Lambda^\varepsilon] S_\varepsilon \|_{L_2(\mathbb{R}^d;\mathbb{C}^m) \to L_2(\mathbb{R}^d;\mathbb{C}^n)}
\leq M.
\eqno(3.9)
$$

\smallskip\noindent\textbf{Lemma 3.5.} \textit{We have}
$$
\| \varepsilon \Lambda^\varepsilon (\Pi_\varepsilon -
S_\varepsilon) b(\mathbf{D}) \mathbf{u}_0 \|_{H^1(\mathbb{R}^d;\mathbb{C}^n)}
\leq \check{C} \varepsilon \|\mathbf{u}_0\|_{H^2(\mathbb{R}^d;\mathbb{C}^n)}.
\eqno(3.10)
$$
\textit{The constant} $\check{C}$ \textit{is defined below in} (3.16) \textit{and depends only on}
$m$, $d$, $\|g\|_{L_\infty}$, $\|g^{-1}\|_{L_\infty}$, $\alpha_0$, $\alpha_1$,
\textit{and the parameters of the lattice} $\Gamma$.

\smallskip\noindent\textbf{Proof.} From (1.3), (3.8), and (3.9) it follows that
$$
\|\varepsilon \Lambda^\varepsilon (\Pi_\varepsilon -
S_\varepsilon) b(\mathbf{D}) \mathbf{u}_0\|_{L_2(\mathbb{R}^d;\mathbb{C}^n)} \leq
2 M \alpha_1^{1/2} \varepsilon \|\mathbf{u}_0 \|_{H^1(\mathbb{R}^d;\mathbb{C}^n)}.
\eqno(3.11)
$$
Consider the derivatives
$$
\begin{aligned}
\frac{\partial}{\partial x_j}(\varepsilon \Lambda^\varepsilon
(\Pi_\varepsilon - S_\varepsilon) b(\mathbf{D}) \mathbf{u}_0)
=& \left( \frac{\partial \Lambda}{\partial x_j}\right)^\varepsilon
(\Pi_\varepsilon - S_\varepsilon) b(\mathbf{D}) \mathbf{u}_0
\\
&+ \varepsilon \Lambda^\varepsilon (\Pi_\varepsilon - S_\varepsilon) b(\mathbf{D})
\partial_j \mathbf{u}_0,\quad j=1,\dots,d.
\end{aligned}
$$
Then
$$
\begin{aligned}
&\sum_{j=1}^d \|\partial_j (\varepsilon \Lambda^\varepsilon (\Pi_\varepsilon - S_\varepsilon)
b(\mathbf{D}) \mathbf{u}_0)\|^2_{L_2(\mathbb{R}^d)}
\\
&\leq 2 \int_{\mathbb{R}^d} |(\mathbf{D} \Lambda)^\varepsilon|^2 |
(\Pi_\varepsilon - S_\varepsilon) b(\mathbf{D}) \mathbf{u}_0|^2\,d\mathbf{x}
\\
&+ 2 \varepsilon^2 \sum_{j=1}^d \int_{\mathbb{R}^d} |\Lambda^\varepsilon (\Pi_\varepsilon - S_\varepsilon)
b(\mathbf{D}) \partial_j \mathbf{u}_0|^2\,d\mathbf{x}.
\end{aligned}
\eqno(3.12)
$$
The second summand in the right-hand side of (3.12) is estimated by using (1.3), (3.8), and (3.9):
$$
2 \varepsilon^2 \sum_{j=1}^d \int_{\mathbb{R}^d} |\Lambda^\varepsilon (\Pi_\varepsilon - S_\varepsilon) b(\mathbf{D}) \partial_j \mathbf{u}_0|^2\,d\mathbf{x}
\leq 8 \varepsilon^2 M^2 \alpha_1 \|\mathbf{u}_0 \|^2_{H^2(\mathbb{R}^d;\mathbb{C}^n)}.
\eqno(3.13)
$$
The first summand in the right-hand side of (3.12) is estimated with the help of Lemma 3.4:
$$
\begin{aligned}
2&\int_{\mathbb{R}^d} |(\mathbf{D} \Lambda)^\varepsilon|^2 |
(\Pi_\varepsilon - S_\varepsilon) b(\mathbf{D})\mathbf{u}_0|^2 \,d\mathbf{x} \leq
2\beta_1 \int_{\mathbb{R}^d}|(\Pi_\varepsilon - S_\varepsilon) b(\mathbf{D})\mathbf{u}_0|^2 \,d\mathbf{x}
\\
&+ 2\beta_2 \varepsilon^2 \sum_{j=1}^d
\int_{\mathbb{R}^d} |\Lambda^\varepsilon|^2 |(\Pi_\varepsilon - S_\varepsilon)
b(\mathbf{D}) \partial_j \mathbf{u}_0|^2\,d\mathbf{x}.
\end{aligned}
\eqno(3.14)
$$
Next, by Propositions 1.4 and 3.1 and relation (1.3), we have
$$
\| (\Pi_\varepsilon - S_\varepsilon) b(\mathbf{D})\mathbf{u}_0
\|_{L_2(\mathbb{R}^d)} \leq \varepsilon (r_0^{-1} + r_1)
\alpha_1^{1/2} \| \mathbf{u}_0 \|_{H^2(\mathbb{R}^d;\mathbb{C}^n)}.
\eqno(3.15)
$$
The second summand in the right-hand side of (3.14) is estimated with the help of (3.13).
Finally, combining (3.12)--(3.15), we obtain
$$
\begin{aligned}
&\sum_{j=1}^d \|\partial_j (\varepsilon \Lambda^\varepsilon (\Pi_\varepsilon - S_\varepsilon)
b(\mathbf{D}) \mathbf{u}_0)\|^2_{L_2(\mathbb{R}^d)}
\\
&\leq \varepsilon^2 \left( 8 M^2 (1+\beta_2) + 2 \beta_1(r_0^{-1} + r_1)^2\right)
\alpha_1 \| \mathbf{u}_0 \|^2_{H^2(\mathbb{R}^d;\mathbb{C}^n)}.
\end{aligned}
$$
Together with (3.11) this implies (3.10) with the constant
$$
\check{C} = \alpha_1^{1/2} \left( M^2 (8\beta_2 + 12) +
2 \beta_1 (r_0^{-1} + r_1)^2 \right)^{1/2}.\ \ \bullet
\eqno(3.16)
$$

\smallskip
Now it is easy to complete the \textbf{proof of Theorem 3.3}.
By (1.3) and (1.10), we obtain the following lower estimate for
the symbol of the effective operator:
$$
b(\boldsymbol{\xi})^* g^0 b(\boldsymbol{\xi}) \geq c_0 |\boldsymbol{\xi}|^2 \boldsymbol{1}_n,
\ \ \boldsymbol{\xi} \in \mathbb{R}^d,\ \ c_0 = \alpha_0 \|g^{-1}\|^{-1}_{L_\infty}.
\eqno(3.17)
$$
Using the Fourier transformation and (3.17), we estimate the norm of the function
$\mathbf{u}_0 = (\mathcal{A}^0 +I)^{-1} \mathbf{F}$ in $H^2(\mathbb{R}^d;\mathbb{C}^n)$:
$$
\begin{aligned}
&\| \mathbf{u}_0 \|^2_{H^2(\mathbb{R}^d;\mathbb{C}^n)} = \int_{\mathbb{R}^d} (1+ |\boldsymbol{\xi}|^2)^2
\left| (b(\boldsymbol{\xi})^* g^0 b(\boldsymbol{\xi}) + \boldsymbol{1}_n)^{-1}
\widehat{\mathbf{F}}(\boldsymbol{\xi})\right|^2 \,d\boldsymbol{\xi}
\\
&\leq
\int_{\mathbb{R}^d} (1+ |\boldsymbol{\xi}|^2)^2
 (c_0 |\boldsymbol{\xi}|^2 +1)^{-2}|\widehat{\mathbf{F}}(\boldsymbol{\xi})|^2 \,d\boldsymbol{\xi}
 \leq (1+ c_0^{-1})^2 \| \mathbf{F} \|^2_{L_2(\mathbb{R}^d;\mathbb{C}^n)}.
\end{aligned}
$$
Combining this with (1.15), (3.5), and (3.10), we obtain
$$
\| \mathbf{v}_\varepsilon - \widetilde{\mathbf{v}}_\varepsilon \|_{H^1(\mathbb{R}^d;\mathbb{C}^n)} \leq \check{C} \varepsilon \| \mathbf{u}_0 \|_{H^2(\mathbb{R}^d;\mathbb{C}^n)}
\leq (1+ c_0^{-1}) \check{C} \varepsilon \| \mathbf{F} \|_{L_2(\mathbb{R}^d;\mathbb{C}^n)}.
\eqno(3.18)
$$
Relations (1.16) and (3.18) imply (3.6) with
$\widetilde{C}_2 = C_2 + (1+ c_0^{-1}) \check{C}$. $\ \bullet$

\section*{\S 4. Homogenization of the Dirichlet problem in a bounded domain:
preliminaries}

\smallskip\noindent\textbf{4.1. Statement of the problem.}
Let $\mathcal{O} \subset \mathbb{R}^{d}$ be a bounded domain of class $C^2$.
In $L_2 (\mathcal{O}; \mathbb{C}^{n})$, we consider the operator
$\mathcal{A}_{D,\varepsilon}$ formally given by the differential expression
$b (\mathbf{D})^* g^{\varepsilon} (\mathbf{x}) b (\mathbf{D})$
with the Dirichlet condition on $\partial \mathcal{O}$.
Precisely, $\mathcal{A}_{D,\varepsilon}$ is the selfadjoint operator in
$L_2(\mathcal{O};\mathbb{C}^n)$ generated by the quadratic form
$$
a_{D,\varepsilon}[\mathbf{u},\mathbf{u}] = \int_{\mathcal{O}}
\left\langle g^{\varepsilon} (\mathbf{x}) b (\mathbf{D}) \mathbf{u},
b(\mathbf{D}) \mathbf{u} \right\rangle \, d \mathbf{x},
\quad \mathbf{u} \in H^1_0 (\mathcal{O}; \mathbb{C}^{n}).
$$
This form is closed and positive definite. Indeed,
let us extend $\mathbf{u}$ by zero to $\mathbb{R}^d \setminus \mathcal{O}$.
Then $\mathbf{u} \in H^1(\mathbb{R}^d;\mathbb{C}^n)$. Applying (1.4), we obtain
$$
c_0 \intop_{\mathcal{O}} |\mathbf{D} \mathbf{u}|^2\, d\mathbf{x}
\leq a_{D,\varepsilon}[\mathbf{u}, \mathbf{u}]
\leq c_1 \intop_{\mathcal{O}} |\mathbf{D} \mathbf{u}|^2\, d\mathbf{x},
\ \ \mathbf{u} \in H^1_0(\mathcal{O};\mathbb{C}^n).
\eqno(4.1)
$$
It remains to note that the functional $\|\mathbf{D} \mathbf{u}\|_{L_2(\mathcal{O})}$
determines the norm in $H^1_0(\mathcal{O};\mathbb{C}^n)$ equivalent to the standard one.

\textit{Our goal} is to find approximation for small $\varepsilon$
for the operator $\mathcal{A}_{D,\varepsilon}^{-1}$
in the norm of operators acting from $L_2(\mathcal{O};\mathbb{C}^n)$ to $H^1(\mathcal{O};\mathbb{C}^n)$.
In terms of solutions, we are interested in the behavior of the generalized solution
$\mathbf{u}_\varepsilon \in H^1_0(\mathcal{O};\mathbb{C}^n)$ of the Dirichlet problem
$$
b (\mathbf{D})^* g^{\varepsilon} (\mathbf{x}) b (\mathbf{D}) \mathbf{u}_{\varepsilon}(\mathbf{x})
= \mathbf{F}(\mathbf{x}),\ \ \mathbf{x} \in \mathcal{O};
\quad \mathbf{u}_\varepsilon\vert_{\partial \mathcal{O}}=0,
\eqno(4.2)
$$
where $\mathbf{F} \in L_2 (\mathcal{O}; \mathbb{C}^{n})$.
Then $\mathbf{u}_\varepsilon = \mathcal{A}_{D,\varepsilon}^{-1} \mathbf{F}$.

\smallskip\noindent\textbf{4.2. The energy inequality.}
Now, we consider the problem (4.2) with the right-hand side of class $H^{-1}(\mathcal{O};\mathbb{C}^n)$
and prove the energy inequality.
Recall that $H^{-1}(\mathcal{O};\mathbb{C}^n)$ is defined as the space dual to
$H^1_0(\mathcal{O};\mathbb{C}^n)$ with respect to the $L_2(\mathcal{O};\mathbb{C}^n)$-coupling.
If $\mathbf{f} \in H^{-1}(\mathcal{O};\mathbb{C}^n)$ and
$\boldsymbol{\eta} \in H^1_0 (\mathcal{O}; \mathbb{C}^{n})$, the symbol
$\int_\mathcal{O} \langle \mathbf{f}, \boldsymbol{\eta} \rangle\,d\mathbf{x}$
stands for the value of the functional $\mathbf{f}$ on the element $\boldsymbol{\eta}$.
Herewith,
$$
\left| \int_\mathcal{O} \langle \mathbf{f}, \boldsymbol{\eta} \rangle\,d\mathbf{x}\right|
\le \|\mathbf{f}\|_{H^{-1}(\mathcal{O};\mathbb{C}^n)} \|\boldsymbol{\eta}\|_{H^1(\mathcal{O};\mathbb{C}^n)}.
\eqno(4.3)
$$

\smallskip\noindent\textbf{Lemma 4.1.} \textit{Let} $\mathbf{f} \in H^{-1}(\mathcal{O};\mathbb{C}^n)$,
\textit{and let} $\mathbf{z}_\varepsilon \in H^1_0(\mathcal{O};\mathbb{C}^n)$
\textit{be the generalized solution of the Dirichlet problem}
$$
b (\mathbf{D})^* g^{\varepsilon} (\mathbf{x}) b (\mathbf{D}) \mathbf{z}_{\varepsilon}(\mathbf{x})
= \mathbf{f}(\mathbf{x}),\ \ \mathbf{x} \in \mathcal{O};
\quad \mathbf{z}_\varepsilon\vert_{\partial \mathcal{O}}=0.
$$
\textit{In other words}, $\mathbf{z}_\varepsilon$ \textit{satisfies the identity}
$$
\int_{\mathcal{O}} \langle g^{\varepsilon}(\mathbf{x}) b (\mathbf{D}) \mathbf{z}_{\varepsilon},
b (\mathbf{D}) \boldsymbol{\eta} \rangle\,d\mathbf{x} = \int_\mathcal{O} \langle \mathbf{f},
\boldsymbol{\eta} \rangle \,d\mathbf{x}, \quad \forall
\ \boldsymbol{\eta} \in H^1_0 (\mathcal{O}; \mathbb{C}^{n}).
\eqno(4.4)
$$
\textit{Then the following estimate called the "energy inequality"\ is true}:
$$
\| \mathbf{z}_{\varepsilon} \|_{H^1 (\mathcal{O}; \mathbb{C}^{n})} \leq
\widehat{C} \| \mathbf{f} \|_{H^{-1} (\mathcal{O}; \mathbb{C}^{n})}.
\eqno(4.5)
$$
\textit{Here} $\widehat{C}=(1 + (\text{diam}\,{\mathcal{O}})^2) \alpha_0^{-1} \| g^{-1} \|_{L_\infty}$.

\smallskip\noindent\textbf{Proof.}
By the lower estimate (4.1), we have
$$
\| \mathbf{D} \mathbf{z}_{\varepsilon} \|^2_{L_2 (\mathcal{O})}
\leq c_0^{-1} \left( g^{\varepsilon} b (\mathbf{D}) \mathbf{z}_{\varepsilon}, b(\mathbf{D}) \mathbf{z}_{\varepsilon}\right)_{L_2 (\mathcal{O})}.
 \eqno(4.6)
 $$
Next, from (4.3) and (4.4) with $\boldsymbol{\eta}= \mathbf{z}_\varepsilon$
it follows that
$$
\left( g^{\varepsilon} b (\mathbf{D}) \mathbf{z}_{\varepsilon}, b(\mathbf{D}) \mathbf{z}_{\varepsilon}
\right)_{L_2 (\mathcal{O})} =  \int_\mathcal{O} \langle \mathbf{f}, \mathbf{z}_{\varepsilon}
\rangle\,d\mathbf{x} \leq \| \mathbf{f} \|_{H^{-1} (\mathcal{O})}
\| \mathbf{z}_{\varepsilon} \|_{H^1 (\mathcal{O})}.
\eqno(4.7)
$$
By the Friedrichs inequality,
$$
\| \mathbf{z}_{\varepsilon} \|_{L_2 (\mathcal{O})} \leq (\text{diam}\,{\mathcal{O}})
\| \mathbf{D} \mathbf{z}_{\varepsilon} \|_{L_2 (\mathcal{O})}.
\eqno(4.8)
$$

Finally, combining relations (4.6)--(4.8), we obtain
$$
\begin{aligned}
\| \mathbf{z}_{\varepsilon} \|^2_{H^1 (\mathcal{O})} &\leq
(1 + (\text{diam}\,\mathcal{O})^2) \| \mathbf{D} \mathbf{z}_{\varepsilon} \|^2_{L_2 (\mathcal{O})}
\\
&\leq (1 + (\text{diam}\,{\mathcal{O}})^2) c_0^{-1}
\| \mathbf{f} \|_{H^{-1} (\mathcal{O})} \| \mathbf{z}_{\varepsilon} \|_{H^1 (\mathcal{O})}.
\end{aligned}
$$
This implies (4.5). $\ \bullet$

\smallskip

Roughening the result of Lemma 4.1, we arrive at the following corollary.

\smallskip\noindent\textbf{Corollary 4.2.} \textit{The operator}
$\mathcal{A}_{D,\varepsilon}^{-1}$ \textit{is continuous from} $L_2(\mathcal{O};\mathbb{C}^n)$
\textit{to} $H^1_0(\mathcal{O};\mathbb{C}^n)$, \textit{and}
$$
\|\mathcal{A}_{D,\varepsilon}^{-1}\|_{L_2(\mathcal{O};\mathbb{C}^n) \to H^1(\mathcal{O};\mathbb{C}^n)}
\leq \widehat{C}.
$$

\smallskip
In what follows, we shall need the next statement which is proved with the help of Lemma 4.1.

\smallskip\noindent\textbf{Lemma 4.3.} \textit{Let}
$\boldsymbol{\psi} \in H^{1}(\mathcal{O};\mathbb{C}^n)$,
\textit{and let} ${\mathbf{r}}_\varepsilon \in H^1(\mathcal{O};\mathbb{C}^n)$
\textit{be the generalized solution of the problem}
$$
b (\mathbf{D})^* g^{\varepsilon} (\mathbf{x}) b (\mathbf{D}) \mathbf{r}_{\varepsilon}(\mathbf{x})
= 0,\ \ \mathbf{x} \in \mathcal{O};
\quad \mathbf{r}_\varepsilon\vert_{\partial \mathcal{O}}= \boldsymbol{\psi} \vert_{\partial \mathcal{O}}.
\eqno(4.9)
$$
\textit{Then}
$$
\| \mathbf{r}_{\varepsilon} \|_{H^1 (\mathcal{O}; \mathbb{C}^{n})} \leq
\gamma_0 \| \boldsymbol{\psi} \|_{H^{1} (\mathcal{O}; \mathbb{C}^{n})},
\ \ \gamma_0 = 1+ \widehat{C} d^{1/2} \alpha_1 \|g\|_{L_\infty}.
\eqno(4.10)
$$

\smallskip\noindent\textbf{Proof.}
By (4.9), the function $\mathbf{r}_{\varepsilon} - {\boldsymbol{\psi}}$ is the solution of the
Dirichlet problem
$$
\mathcal{A}_{\varepsilon} (\mathbf{r}_{\varepsilon} -{\boldsymbol{\psi}})
= - \mathcal{A}_{\varepsilon} {\boldsymbol{\psi}} \ \ \text{in}\ \mathcal{O}; \quad
(\mathbf{r}_{\varepsilon} -{\boldsymbol{\psi}}) |_{\partial \mathcal{O}} = 0.
\eqno(4.11)
$$
Here the right-hand side in the equation belongs to $H^{-1}(\mathcal{O};\mathbb{C}^n)$, and
$$
\begin{aligned}
&\| \mathcal{A}_{\varepsilon} {\boldsymbol{\psi}} \|_{H^{-1} (\mathcal{O})} =
\sup_{0\ne \boldsymbol{\varphi} \in H^1_0(\mathcal{O}; \mathbb{C}^n)}
\frac{\left| \left( g^{\varepsilon} b (\mathbf{D}) {\boldsymbol{\psi}},
b (\mathbf{D}) \boldsymbol{\varphi} \right)_{L_2 (\mathcal{O})} \right|}
{\| \boldsymbol{\varphi} \|_{H^1 (\mathcal{O})}}
\\
&\leq \alpha_1^{1/2} \| g \|_{L_\infty}
\| b (\mathbf{D}) {\boldsymbol{\psi}} \|_{L_2 (\mathcal{O})}.
\end{aligned}
\eqno(4.12)
$$
We have taken into account that
$\|b(\mathbf{D}) \boldsymbol{\varphi}\|_{L_2(\mathcal{O})} \leq
\alpha_1^{1/2} \|\mathbf{D} \boldsymbol{\varphi}\|_{L_2(\mathcal{O})}$
which can be checked as follows.
Extend $\boldsymbol{\varphi} \in H^1_0(\mathcal{O};\mathbb{C}^n)$ by zero to
$\mathbb{R}^d \setminus \mathcal{O}$, keeping the same notation $\boldsymbol{\varphi}$.
Then $\boldsymbol{\varphi} \in H^1(\mathbb{R}^d;\mathbb{C}^n)$.
Using the Fourier transformation and the upper inequality (1.3), we obtain
$$
\begin{aligned}
&\|b(\mathbf{D}) \boldsymbol{\varphi}\|_{L_2(\mathcal{O})}^2 =
\|b(\mathbf{D}) \boldsymbol{\varphi}\|_{L_2(\mathbb{R}^d)}^2 =
\int_{\mathbb{R}^d} |b(\boldsymbol{\xi})  \widehat{\boldsymbol{\varphi}}(\boldsymbol{\xi})|^2\,
d\boldsymbol{\xi}
\\
&\leq \alpha_1 \int_{\mathbb{R}^d} |\boldsymbol{\xi}|^2
|\widehat{\boldsymbol{\varphi}}(\boldsymbol{\xi})|^2\, d\boldsymbol{\xi}=
\alpha_1 \|\mathbf{D} \boldsymbol{\varphi}\|_{L_2(\mathbb{R}^d)}^2 =
\alpha_1 \|\mathbf{D} \boldsymbol{\varphi}\|_{L_2(\mathcal{O})}^2.
\end{aligned}
\eqno(4.13)
$$

Next, by (1.2) and (2.6),
$$
\| b (\mathbf{D}) {\boldsymbol{\psi}} \|_{L_2 (\mathcal{O})}
\leq \alpha_1^{1/2} \sum_{l=1}^d \| D_l {\boldsymbol{\psi}} \|_{L_2 (\mathcal{O})}
\leq \alpha_1^{1/2} d^{1/2} \| {\boldsymbol{\psi}} \|_{H^1 (\mathcal{O})}.
\eqno(4.14)
$$
From (4.12) and (4.14) it follows that
$$
\| \mathcal{A}_{\varepsilon} {\boldsymbol{\psi}}
\|_{H^{-1} (\mathcal{O})}
\leq \alpha_1 d^{1/2} \|g\|_{L_\infty} \| {\boldsymbol{\psi}} \|_{H^1 (\mathcal{O})}.
\eqno(4.15)
$$

Applying Lemma 4.1 to the problem (4.11), we obtain
$$
\| \mathbf{r}_{\varepsilon} - {\boldsymbol{\psi}}
\|_{H^1 (\mathcal{O})} \leq \widehat{C}
\| \mathcal{A}_{\varepsilon} {\boldsymbol{\psi}}
\|_{H^{-1} (\mathcal{O})}.
\eqno(4.16)
$$
Now, (4.15) and (4.16) imply (4.10).
$\ \bullet$

\smallskip\noindent\textbf{Remark 4.4.} The statements of Lemma 4.1 and Corollary 4.2 remain true in any
bounded domain $\mathcal{O} \subset \mathbb{R}^d$ (without the assumption that $\partial \mathcal{O} \in C^2$).
The same is true for Lemma 4.3 if the problem (4.9) is understood as the identity
$$
\int_\mathcal{O} \langle g^\varepsilon(\mathbf{x}) b(\mathbf{D})\mathbf{r}_{\varepsilon},b(\mathbf{D})
\boldsymbol{\eta}\rangle \,d\mathbf{x} =0,\quad \forall \boldsymbol{\eta}
\in H^1_0(\mathcal{O};\mathbb{C}^n),
$$
and relation
$\mathbf{r}_{\varepsilon} -{\boldsymbol{\psi}} \in H^1_0(\mathcal{O};\mathbb{C}^n)$.

\smallskip\noindent\textbf{4.3. The "homogenized"\ problem.}
In $L_2(\mathcal{O};\mathbb{C}^n)$, we consider the selfadjoint operator
$\mathcal{A}_D^0$ generated by the quadratic form
$$
\int_{\mathcal{O}} \left\langle g^0 b (\mathbf{D}) \mathbf{u},
b(\mathbf{D}) \mathbf{u} \right\rangle \,d \mathbf{x},
 \quad \mathbf{u} \in H^1_0 (\mathcal{O}; \mathbb{C}^{n}).
$$
Here $g^0$ is the effective matrix defined by (1.6).
Applying Corollary 4.2 with $g^\varepsilon$ replaced by $g^0$ and taking (1.10) into account,
we see that the operator $(\mathcal{A}^0_D)^{-1}$ is continuous from
$L_2(\mathcal{O};\mathbb{C}^n)$ to $H^1_0(\mathcal{O};\mathbb{C}^n)$, and
$$
\| (\mathcal{A}^0_D)^{-1}\|_{L_2(\mathcal{O};\mathbb{C}^n) \to H^1(\mathcal{O};\mathbb{C}^n)}
\leq \widehat{C},
\eqno(4.17)
$$
where the constant $\widehat{C}$ is defined in Lemma 4.1.
Note that this fact is valid in any bounded domain
$\mathcal{O} \subset \mathbb{R}^d$ (without the assumption that $\partial \mathcal{O} \in C^2$).

Let $\mathbf{u}_0 \in H^1_0(\mathcal{O};\mathbb{C}^n)$ be the generalized solution of the Dirichlet
problem
$$
b(\mathbf{D})^* g^0 b(\mathbf{D}) \mathbf{u}_0(\mathbf{x}) = \mathbf{F}(\mathbf{x}),
\ \ \mathbf{x} \in \mathcal{O};\ \ \mathbf{u}_0\vert_{\partial \mathcal{O}} =0,
\eqno(4.18)
$$
where $\mathbf{F} \in L_2(\mathcal{O};\mathbb{C}^n)$. Then
$\mathbf{u}_0 = (\mathcal{A}^0_D)^{-1} \mathbf{F}$.

Since $\partial \mathcal{O} \in C^2$, for the solution $\mathbf{u}_0$ of the problem (4.18)
we have
$
\mathbf{u}_0 \in H^1_0 (\mathcal{O}; \mathbb{C}^{n}) \cap H^2 (\mathcal{O}; \mathbb{C}^{n}),
$
and
$$
\| \mathbf{u}_0 \|_{H^2 (\mathcal{O}; \mathbb{C}^{n})} \leq \widehat{c}
\| \mathbf{F} \|_{L_2 (\mathcal{O}; \mathbb{C}^{n})}.
\eqno(4.19)
$$
Here the constant $\widehat{c}$ depends only on $\alpha_0$, $\alpha_1$,
$\|g\|_{L_\infty}$, $\|g^{-1}\|_{L_\infty}$, and the domain $\mathcal{O}$.
To justify these properties, it suffices to note that the operator $b(\mathbf{D})^* g^0 b(\mathbf{D})$ is
a \textit{strongly elliptic} matrix DO
and to apply the "additional smoothness"\ theorems for solutions of strongly elliptic systems
(see, e.~g., [McL, Chapter 4]).

It follows that the operator $\mathcal{A}_D^0$ is given by the differential expression
$b(\mathbf{D})^* g^0 b(\mathbf{D})$ on the domain
$H^1_0 (\mathcal{O}; \mathbb{C}^{n}) \cap H^2 (\mathcal{O}; \mathbb{C}^{n})$,
and that the inverse operator satisfies the estimate
$$
\| (\mathcal{A}^0_D)^{-1}\|_{L_2(\mathcal{O};\mathbb{C}^n) \to H^2(\mathcal{O};\mathbb{C}^n)}
\leq \widehat{c}.
\eqno(4.20)
$$

Below we shall see that
the solution $\mathbf{u}_\varepsilon$ of the problem (4.2) converges in
$L_2(\mathcal{O};\mathbb{C}^n)$ to the solution $\mathbf{u}_0$ of the "homogenized"\
problem (4.18), as $\varepsilon \to 0$. Our \textit{main goal} is to
find approximation for $\mathbf{u}_\varepsilon$
in the norm of $H^1(\mathcal{O};\mathbb{C}^n)$; for this, it is necessary to take
the first order corrector into account.

\section*{\S 5. Auxiliary statements}
In this section, we prove several auxiliary statements needed for further considerations.

\smallskip\noindent\textbf{Lemma 5.1.}
\textit{Let} $\mathcal{O} \subset \mathbb{R}^d$ \textit{be a bounded domain of class} $C^1$.
\textit{Denote}
$B_{\varepsilon} = \left\{ \mathbf{x} \in \mathcal{O}: \textrm{\upshape dist\itshape}
\,\{ \mathbf{x}, \partial \mathcal{O} \} < \varepsilon \right\}$.
\textit{Then there exists a number} $\varepsilon_0 \in (0,1]$
\textit{depending on the domain} $\mathcal{O}$ \textit{such that for any}
$u \in H^1 (\mathcal{O})$ \textit{we have}
$$
\int_{B_{\varepsilon}} |u|^{2} d \mathbf{x}
\leq \beta \varepsilon \| u \|_{H^1 (\mathcal{O})} \| u \|_{L_2 (\mathcal{O})},
\quad 0< \varepsilon \leq \varepsilon_0.
\eqno(5.1)
$$
\textit{The constant} $\beta=\beta(\mathcal{O})$ \textit{depends only on the domain} $\mathcal{O}$.

\smallskip\noindent\textbf{Proof.}
Let us start with a model problem in the semiball
${\mathcal{D}}_0 = \left\{ \mathbf{x} \in \mathbb{R}^d: \ |\mathbf{x}| < 1, \ x_d > 0 \right\}$.
For points $\mathbf{x}\in \mathbb{R}^d$ we write $\mathbf{x} = (\mathbf{x}',x_d)$,
where $\mathbf{x'} = (x_1, \ldots, x_{d-1})$.
Introduce the following notation:
$$
\begin{aligned}
{\mathcal{D}}_t &= \{\mathbf{x} \in \mathbb{R}^d:\ |\mathbf{x}|<1,\ x_d >t\},\ \
\Sigma_t = \{ \mathbf{x}\in \partial {\mathcal{D}}_t: \ x_d=t \},\ 0\le t \le \varepsilon;
\\
\Upsilon_\varepsilon &=\{\mathbf{x} \in \mathbb{R}^d:\ |\mathbf{x}|<1,\ 0< x_d < \varepsilon \},
\ \ \Sigma = \{ \mathbf{x} \in \partial {\mathcal{D}}_0: \ |\mathbf{x}|=1\}.
\end{aligned}
$$
Assume that $u \in H^1 ({\mathcal{D}}_0)$ and $u = 0$ on $\Sigma$.
Let $0\leq t\leq \varepsilon$. Using the Green formula in the domain ${\mathcal{D}}_t$, we have:
$$
\int_{{\mathcal{D}}_t} \frac{\partial u}{\partial x_d} \, \overline{u} \, d \mathbf{x'} dx_d
= - \int_{\Sigma_t} |u|^{2} \, d \mathbf{x'}
- \int_{{\mathcal{D}}_t} u \, \frac{\partial \overline{u}}{\partial x_d} \, d \mathbf{x'} dx_d.
$$
Hence,
$$
\begin{aligned}
&\int_{\Sigma_t} |u(\mathbf{x}',t)|^{2} \, d \mathbf{x'} \leq
\int_{{\mathcal{D}}_t} 2 \left| \frac{\partial u}{\partial x_d} \right| \,|u| \, d \mathbf{x}
\\
&\leq 2 \left( \int_{{\mathcal{D}}_0} \left| \frac{\partial u}{\partial x_d} \right|^2
\, d \mathbf{x} \right)^{1/2} \left( \int_{{\mathcal{D}}_0} |u|^{2} \,d \mathbf{x} \right)^{1/2}.
\end{aligned}
$$
Integrating over $t\in (0,\varepsilon)$, we obtain
$$
\int_{\Upsilon_\varepsilon} |u|^{2}\,d\mathbf{x} \le 2 \varepsilon
\left( \int_{{\mathcal{D}}_0} \left| \frac{\partial u}{\partial x_d} \right|^2
\, d \mathbf{x} \right)^{1/2} \left( \int_{{\mathcal{D}}_0} |u|^{2} \,d \mathbf{x} \right)^{1/2}.
$$

Estimate (5.1) in the case of a bounded domain $\mathcal{O}$ of class $C^1$ is deduced
from here in a standard way with the help of local maps,
diffeomorphisms rectifying the boundary, and the partition of unity.
Herewith, we take into account that the space $H^1$ is invariant with respect to
diffeomorphisms of class $C^1$.
The number $\varepsilon_0$ must be such that the set $B_{\varepsilon_0}$
can be covered by a finite number of open sets admitting diffeomorphisms
rectifying the boundary.
Thus, the number $\varepsilon_0$ depends only on the domain $\mathcal{O}$.
$\bullet$

\smallskip

Next statement is a direct consequence of Lemma 5.1.

\smallskip\noindent\textbf{Lemma 5.2.}
\textit{Let} $\mathcal{O} \subset \mathbb{R}^d$ \textit{be a bounded domain of class} $C^1$.
\textit{Denote} $\left(\partial \mathcal{O} \right)_{\varepsilon} =
\left\{ \mathbf{x} \in \mathbb{R}^d: \textrm{dist} \,
\{ \mathbf{x}, \partial \mathcal{O} \} < \varepsilon \right\}$.
\textit{Let} $\varepsilon_1 \in (0,1]$ \textit{be such that the set} $(\partial \mathcal{O})_{\varepsilon_1}$
\textit{can be covered by a finite number of open sets admitting diffeomorphisms of class} $C^1$
\textit{rectifying the boundary} $\partial \mathcal{O}$.
\textit{Then for any} $u \in H^1 (\mathbb{R}^d)$ \textit{we have}
$$
\int_{\left( \partial \mathcal{O} \right)_{\varepsilon}} |u|^{2} \,d \mathbf{x}
\leq \beta^0 \varepsilon \| u \|_{H^1 (\mathbb{R}^d)} \| u \|_{L_2 (\mathbb{R}^d)},
\ \ 0< \varepsilon \leq \varepsilon_1.
\eqno(5.2)
$$
\textit{The constant} $\beta^0=\beta^0(\mathcal{O})$ \textit{depends only on the domain} $\mathcal{O}$.

\smallskip
\noindent\textbf{Proof.} We apply Lemma 5.1 in the domain $\mathcal{O}$
and in the domain $\mathcal{B} \setminus \overline{\mathcal{O}}$, where
$\mathcal{B}$ is some open ball containing
$\overline{\mathcal{O}} \cup \overline{(\partial \mathcal{O})_{\varepsilon_1}}$.
Then (5.2) is true with
$\beta^0 = \max \{ \beta(\mathcal{O}), \beta(\mathcal{B} \setminus \overline{\mathcal{O}})\}$.
$\ \bullet$

The following statement is similar to Lemma 2.6 from [ZhPas].

\smallskip\noindent\textbf{Lemma 5.3.} \textit{Let} $S_\varepsilon$ \textit{be the operator} (3.1).
\textit{Suppose that the domain} $\mathcal{O}$ \textit{and the number}
$\varepsilon_1$ \textit{satisfy the assumptions of Lemma} 5.2.
\textit{Assume that} $f(\mathbf{x})$ \textit{is a} $\Gamma$-\textit{periodic function in} $\mathbb{R}^d$
\textit{such that} $f \in L_2(\Omega)$.
\textit{Then for any} $\mathbf{u} \in H^1(\mathbb{R}^d;\mathbb{C}^m)$ \textit{we have}
$$
\begin{aligned}
\int_{\left( \partial \mathcal{O} \right)_{\varepsilon}} |f^\varepsilon(\mathbf{x})|^2
|(S_\varepsilon \mathbf{u})(\mathbf{x})|^{2} \,d \mathbf{x}
\leq \beta_* \varepsilon |\Omega|^{-1}\|f\|^2_{L_2(\Omega)} \| \mathbf{u} \|_{H^1 (\mathbb{R}^d;\mathbb{C}^m)} \| \mathbf{u} \|_{L_2 (\mathbb{R}^d;\mathbb{C}^m)},
\\
0< \varepsilon \leq \varepsilon_2,
\end{aligned}
\eqno(5.3)
$$
\textit{where} $\varepsilon_2 = \varepsilon_1 (1+ r_1)^{-1}$,
$\beta_*= \beta^0 (1+ r_1)$, $2r_1 =\text{diam}\,\Omega$.

\smallskip\noindent\textbf{Proof.} From (3.1), by the Cauchy inequality and the change of variables,
we obtain
$$
\begin{aligned}
&\int_{\left( \partial \mathcal{O} \right)_{\varepsilon}} |f^\varepsilon(\mathbf{x})|^2
|(S_\varepsilon \mathbf{u})(\mathbf{x})|^{2} \,d \mathbf{x}
\leq |\Omega|^{-1} \int_{\left( \partial \mathcal{O} \right)_{\varepsilon}}
d\mathbf{x} \,|f(\varepsilon^{-1}\mathbf{x})|^2
\int_\Omega |\mathbf{u}(\mathbf{x} -\varepsilon \mathbf{z})|^2\,d\mathbf{z}
\\
&\leq |\Omega|^{-1} \int_{\left( \partial \mathcal{O} \right)_{\widetilde{\varepsilon}}} d\mathbf{y}
\int_\Omega d\mathbf{z}\,
|f(\varepsilon^{-1}\mathbf{y} + \mathbf{z})|^2 |\mathbf{u}(\mathbf{y})|^2
\\
&\leq
|\Omega|^{-1} \|f\|^2_{L_2(\Omega)}\int_{\left( \partial \mathcal{O} \right)_{\widetilde{\varepsilon}}}
|\mathbf{u}(\mathbf{y})|^2\,d\mathbf{y}.
\end{aligned}
$$
Here $\widetilde{\varepsilon} = \varepsilon (1+ r_1)$.
Applying Lemma 5.2, we arrive at (5.3). $\ \bullet$

\section*{\S 6. Results in the case of bounded $\Lambda$}

\smallskip\noindent\textbf{6.1.}
We start with the case where Condition 1.9 is satisfied. Denote
$$
K_D^0(\varepsilon) = [\Lambda^{\varepsilon}] b (\mathbf{D}) (\mathcal{A}^0_D)^{-1}.
\eqno(6.1)
$$
By (4.20), the operator $b(\mathbf{D}) (\mathcal{A}^0_D)^{-1}$ is a continuous mapping of
$L_2(\mathcal{O};\mathbb{C}^n)$ into $H^1(\mathcal{O};\mathbb{C}^m)$.
Under Condition 1.9 the operator $[\Lambda^\varepsilon]$
of multiplication by the matrix-valued function $\Lambda^\varepsilon(\mathbf{x})$
is continuous from $H^1(\mathcal{O};\mathbb{C}^m)$ to $H^1(\mathcal{O};\mathbb{C}^n)$.
This easily follows from Corollary 2.4. Consequently, the operator (6.1)
is continuous from $L_2(\mathcal{O};\mathbb{C}^n)$ to $H^1(\mathcal{O};\mathbb{C}^n)$.

Let $\mathbf{u}_\varepsilon$ be the solution of the problem (4.2),
and let $\mathbf{u}_0$ be the solution of the problem (4.18).
The "first order approximation"\ of $\mathbf{u}_\varepsilon$ is given by
$$
\mathbf{\check{v}}_{\varepsilon}= \mathbf{u}_0 + \varepsilon \Lambda^\varepsilon
b (\mathbf{D}) \mathbf{u}_0 = (\mathcal{A}_D^0)^{-1}\mathbf{F} + \varepsilon K_D^0(\varepsilon) \mathbf{F}.
\eqno(6.2)
$$

The following theorem is \textit{our main result} in the case where $\Lambda \in L_\infty$.

\smallskip\noindent\textbf{Theorem 6.1.} \textit{Suppose that} $\mathcal{O}\subset \mathbb{R}^d$
\textit{is a bounded domain of class} $C^2$. \textit{Let} $g(\mathbf{x})$
\textit{and} $b(\mathbf{D})$ \textit{satisfy the assumptions of Subsection} 1.2.
\textit{Let} $\mathbf{u}_{\varepsilon}$ \textit{be the solution of the problem} (4.2), \textit{and let}
$\mathbf{u}_0$ \textit{be the solution of the problem} (4.18) \textit{with}
$\mathbf{F} \in L_2(\mathcal{O};\mathbb{C}^n)$.
\textit{Suppose that} $\Lambda (\mathbf{x})$
\textit{is the} $\Gamma$-\textit{periodic solution of the problem} (1.5)
\textit{and Condition} 1.9 \textit{is satisfied}.
\textit{Let} $\mathbf{\check{v}}_{\varepsilon}$
\textit{be the function defined by} (6.2). \textit{Then there exists a number} $\varepsilon_1 \in (0,1]$
\textit{depending on the domain} $\mathcal{O}$ \textit{such that we have}
$$
\| \mathbf{u}_{\varepsilon} - \mathbf{\check{v}}_{\varepsilon}
\|_{H^1 (\mathcal{O}; \mathbb{C}^{n})} \leq C_0 \varepsilon^{1/2}
\| \mathbf{F} \|_{L_2 (\mathcal{O}; \mathbb{C}^{n})},\ \ 0< \varepsilon \leq \varepsilon_1,
\eqno(6.3)
$$
\textit{or, in operator terms,}
$$
\| \mathcal{A}_{D,\varepsilon}^{-1} - (\mathcal{A}^0_D)^{-1} -
\varepsilon K_D^0 (\varepsilon) \|_{L_2 (\mathcal{O}; \mathbb{C}^{n}) \to
H^1 (\mathcal{O}; \mathbb{C}^{n})} \leq C_0 {\varepsilon}^{1/2},\ \ 0< \varepsilon \leq \varepsilon_1.
$$
\textit{The flux} $\mathbf{p}_\varepsilon:= g^\varepsilon b(\mathbf{D})\mathbf{u}_\varepsilon$
\textit{admits the following approximation}
$$
\| \mathbf{p}_\varepsilon - \widetilde{g}^\varepsilon b(\mathbf{D})\mathbf{u}_0 \|_{L_2(\mathcal{O};\mathbb{C}^m)}
\leq C_0' \varepsilon^{1/2} \| \mathbf{F} \|_{L_2 (\mathcal{O}; \mathbb{C}^{n})},\ \ 0< \varepsilon \leq \varepsilon_1,
\eqno(6.4)
$$
\textit{where} $\widetilde{g}(\mathbf{x}):= g(\mathbf{x})(b(\mathbf{D})\Lambda(\mathbf{x}) + \mathbf{1}_m)$.
\textit{The constants} $C_0$, $C_0'$ \textit{depend only on}
$m$, $d$, $\alpha_0$, $\alpha_1$, $\| g \|_{L_\infty}$,
$\| g^{-1} \|_{L_\infty}$, \textit{the parameters of the lattice} $\Gamma$, \textit{the norm}
$\| \Lambda \|_{L_\infty}$, \textit{and the domain} $\mathcal{O}$.

\smallskip
Recall that some sufficient conditions under which Condition 1.9 is satisfied
are given above in Proposition 1.11.
In particular, the statements of Theorem 6.1 are true for all operators of the form (1.1)
in dimension $d \leq 2$, and also for the scalar elliptic operator
$\mathcal{A}_\varepsilon = - \text{div}\, g^\varepsilon(\mathbf{x}) \nabla$ in arbitrary dimension,
where $g(\mathbf{x})$ is a matrix with real entries.

Roughening the result of Theorem 6.1, we arrive at the following corollary.

\smallskip\noindent\textbf{Corollary 6.2.} \textit{Under the assumptions of Theorem} 6.1
\textit{we have}
$$
\| \mathbf{u}_{\varepsilon} - \mathbf{u}_0
\|_{L_2 (\mathcal{O}; \mathbb{C}^{n})} \leq \widetilde{C}_0 \varepsilon^{1/2}
\| \mathbf{F} \|_{L_2 (\mathcal{O}; \mathbb{C}^{n})},\ \ 0< \varepsilon \leq \varepsilon_1,
\eqno(6.5)
$$
\textit{or, in operator terms},
$$
\| \mathcal{A}_{D,\varepsilon}^{-1} - (\mathcal{A}^0_D)^{-1}
\|_{L_2 (\mathcal{O}; \mathbb{C}^{n}) \to
L_2(\mathcal{O}; \mathbb{C}^{n})} \leq \widetilde{C}_0 {\varepsilon}^{1/2},
\ \ 0< \varepsilon \leq \varepsilon_1.
$$
\textit{Here} $\widetilde{C}_0 = C_0 + \widehat{C} \alpha_1^{1/2} \|\Lambda\|_{L_\infty}$,
\textit{where} $\widehat{C}$ \textit{is defined in Lemma} 4.1.

\smallskip\noindent\textbf{Proof.} From (6.2) and (6.3) it follows that
$$
\| \mathbf{u}_\varepsilon - \mathbf{u}_0 \|_{L_2(\mathcal{O})} \le C_0 \varepsilon^{1/2} \|\mathbf{F}\|_{L_2(\mathcal{O})}+
\varepsilon \| \Lambda^\varepsilon b(\mathbf{D}) \mathbf{u}_0 \|_{L_2(\mathcal{O})},\ \ 0< \varepsilon \leq \varepsilon_1.
\eqno(6.6)
$$
Under Condition 1.9 we have:
$$
\| \Lambda^\varepsilon b(\mathbf{D}) \mathbf{u}_0 \|_{L_2(\mathcal{O})} \le \|\Lambda\|_{L_\infty} \|b(\mathbf{D}) \mathbf{u}_0\|_{L_2(\mathcal{O})}.
\eqno(6.7)
$$
Similarly to (4.13),
$$
\|b(\mathbf{D}) \mathbf{u}_0\|_{L_2(\mathcal{O})} \leq \alpha_1^{1/2} \|\mathbf{D} \mathbf{u}_0\|_{L_2(\mathcal{O})}.
\eqno(6.8)
$$
Combining (6.7) and (6.8) and taking (4.17) into account, we obtain
$$
\| \Lambda^\varepsilon b(\mathbf{D}) \mathbf{u}_0 \|_{L_2(\mathcal{O})} \leq \alpha_1^{1/2} \|\Lambda\|_{L_\infty}
\|\mathbf{D} \mathbf{u}_0\|_{L_2(\mathcal{O})} \leq \widehat{C}\alpha_1^{1/2} \|\Lambda\|_{L_\infty}
\|\mathbf{F}\|_{L_2(\mathcal{O})}.
$$
Together with (6.6) this implies (6.5). \ $\bullet$

\smallskip
Now we distinguish the special cases. Next statement follows from Theorem 6.1
and Propositions 1.2 and 1.3.

\smallskip\noindent\textbf{Proposition 6.3.} $1^\circ$. \textit{If} $g^0 = \overline{g}$,
\textit{i.~e., relations} (1.8) \textit{are satisfied, then} $\Lambda=0$
\textit{and} $K_D^0(\varepsilon)=0$. \textit{In this case we have}
$$
\| \mathbf{u}_{\varepsilon} - \mathbf{u}_0
\|_{H^1 (\mathcal{O}; \mathbb{C}^{n})} \leq {C}_0 \varepsilon^{1/2}
\| \mathbf{F} \|_{L_2 (\mathcal{O}; \mathbb{C}^{n})},\ \ 0< \varepsilon \leq \varepsilon_1.
$$
$2^\circ$. \textit{If} $g^0 = \underline{g}$,
\textit{i.~e., relations} (1.9) \textit{are satisfied, then} $\widetilde{g}=g^0$.
\textit{In this case we have}
$$
\| \mathbf{p}_\varepsilon - {g}^0 b(\mathbf{D})\mathbf{u}_0 \|_{L_2(\mathcal{O};\mathbb{C}^m)}
\leq C_0' \varepsilon^{1/2} \| \mathbf{F} \|_{L_2 (\mathcal{O}; \mathbb{C}^{n})},
\ \ 0< \varepsilon \leq \varepsilon_1.
$$

\smallskip\noindent\textbf{6.2.} The proof of Theorem 6.1
relies on the results for homogenization problem in $\mathbb{R}^d$ (Theorems 1.6 and 1.10)
and on the tricks suggested in [Zh2, ZhPas]
that allow one to carry over such results to the case of a bounded domain.

Let us fix a linear continuous extension operator
$$
P_\mathcal{O}: H^2(\mathcal{O};\mathbb{C}^n) \to H^2(\mathbb{R}^d;\mathbb{C}^n),
\eqno(6.9)
$$
and put $\widetilde{\mathbf{u}}_0 = P_\mathcal{O} \mathbf{u}_0$. Then
$$
\|\widetilde{\mathbf{u}}_0 \|_{H^2(\mathbb{R}^d;\mathbb{C}^n)}
\leq C_\mathcal{O} \| \mathbf{u}_0 \|_{H^2(\mathcal{O};\mathbb{C}^n)},
\eqno(6.10)
$$
where $C_\mathcal{O}$ is the norm of the operator (6.9). Denote
$$
{\mathbf{v}}^{(1)}_{\varepsilon} (\mathbf{x}) =
\widetilde{\mathbf{u}}_0 (\mathbf{x}) + \varepsilon \Lambda^\varepsilon({\mathbf{x}})
b(\mathbf{D}) \widetilde{\mathbf{u}}_0 (\mathbf{x}).
\eqno(6.11)
$$
Then $\mathbf{\check{v}}_{\varepsilon} = {\mathbf{v}}^{(1)}_{\varepsilon} |_{\mathcal{O}}.$

The following statement is proved with the help of Theorems 1.6 and 1.10.

\smallskip\noindent\textbf{Lemma 6.4.} \textit{Let} $\mathbf{u}_0$ \textit{be the solution of
the problem} (4.18), \textit{and let} $\mathbf{\check{v}}_{\varepsilon}$
\textit{be the function defined by} (6.2). \textit{Then for}
$0< \varepsilon \leq 1$ \textit{we have}
$$
\| \mathcal{A}_{\varepsilon} \mathbf{\check{v}}_{\varepsilon}
- \mathcal{A}^0 \mathbf{u}_0 \|_{H^{-1} (\mathcal{O}; \mathbb{C}^{n})}
\leq  C_4 \varepsilon \| \mathbf{u}_0 \|_{H^2 (\mathcal{O}; \mathbb{C}^{n})}.
\eqno(6.12)
$$
\textit{The constant} $C_4$ \textit{depends only on} $m$, $d$, $\alpha_0$, $\alpha_1$,
$\| g \|_{L_\infty}$, $\| g^{-1} \|_{L_\infty}$,
\textit{the parameters of the lattice} $\Gamma$, \textit{the norm} $\| \Lambda \|_{L_\infty}$,
\textit{and the domain} $\mathcal{O}$.

\smallskip\noindent\textbf{Proof.} The required estimate in the case of a bounded domain
is deduced from the similar inequality in $\mathbb{R}^d$. Let
${\mathbf{v}}^{(1)}_{\varepsilon}$ be defined by (6.11). We check that
$$
\| \mathcal{A}_{\varepsilon} {\mathbf{v}}^{(1)}_{\varepsilon} -
\mathcal{A}^0 \widetilde{\mathbf{u}}_0 \|_{H^{-1} (\mathbb{R}^d; \mathbb{C}^{n})} \leq
\widetilde{C}_4 \varepsilon \| \widetilde{\mathbf{u}}_0
\|_{H^2 (\mathbb{R}^d; \mathbb{C}^{n})}, \ \ 0 < \varepsilon \leq 1.
\eqno(6.13)
$$
Clearly,
$$
\widetilde{\mathbf{F}}:= \mathcal{A}^0 \widetilde{\mathbf{u}}_0 + \widetilde{\mathbf{u}}_0
\in L_2 (\mathbb{R}^d; \mathbb{C}^{n}).
\eqno(6.14)
$$
By the Fourier transformation and (1.3), (1.10), we obtain
$$
\begin{aligned}
&\| \widetilde{\mathbf{F}} \|^2_{L_2 (\mathbb{R}^d)}
= \int_{\mathbb{R}^d} \left| (b (\boldsymbol{\xi})^*
g^0 b (\boldsymbol{\xi}) + \mathbf{1}) \widehat{\mathbf{u}}_0 (\boldsymbol{\xi}) \right|^2
\, d \boldsymbol{\xi}
\\
&\leq \ \int_{\mathbb{R}^d} (\alpha_1 |g^0| |\boldsymbol{\xi}|^2 + 1)^2
| \widehat{\mathbf{u}}_0 (\boldsymbol{\xi}) |^2 \, d \boldsymbol{\xi}
\leq \left( \textrm{max} \left\{ \alpha_1 \| g \|_{L_\infty}, 1 \right\} \right)^2
\| \widetilde{\mathbf{u}}_0 \|^2_{H^2 (\mathbb{R}^d)}.
\end{aligned}
\eqno(6.15)
$$
Here $\widehat{\mathbf{u}}_0 (\boldsymbol{\xi})$ is the Fourier-image of the function
$\widetilde{\mathbf{u}}_0 (\mathbf{x})$.

Let $\mathbf{s}_{\varepsilon}\in H^1(\mathbb{R}^d;\mathbb{C}^n)$ be the generalized solution
of the equation
$$
\mathcal{A}_{\varepsilon} \mathbf{s}_{\varepsilon} + \mathbf{s}_{\varepsilon}
= \widetilde{\mathbf{F}}.
\eqno(6.16)
$$
From Theorems 1.6 and 1.10 it follows that
$$
\| \mathbf{s}_{\varepsilon} - \widetilde{\mathbf{u}}_0
\|_{L_2 (\mathbb{R}^d; \mathbb{C}^{n})}  \leq  C_1 \varepsilon
\| \widetilde{\mathbf{F}} \|_{L_2 (\mathbb{R}^d; \mathbb{C}^{n})}, \ \ 0 < \varepsilon \leq 1,
\eqno(6.17)
$$
$$
\| \mathbf{s}_{\varepsilon} - {\mathbf{v}}^{(1)}_{\varepsilon}
\|_{H^1 (\mathbb{R}^d; \mathbb{C}^{n})} \leq C_3 \varepsilon \| \widetilde{\mathbf{F}}
\|_{L_2 (\mathbb{R}^d; \mathbb{C}^{n})}, \ \ 0 < \varepsilon \leq 1.
\eqno(6.18)
$$

By (6.14) and (6.16),
$$
\mathcal{A}_{\varepsilon} {\mathbf{v}}^{(1)}_{\varepsilon} -
\mathcal{A}^0 \widetilde{\mathbf{u}}_0 =
\mathcal{A}_{\varepsilon} ( {\mathbf{v}}^{(1)}_{\varepsilon} - \mathbf{s}_{\varepsilon} ) +
\mathcal{A}_{\varepsilon} \mathbf{s}_{\varepsilon} - \mathcal{A}^0 \widetilde{\mathbf{u}}_0
= \mathcal{A}_{\varepsilon} ( {\mathbf{v}}^{(1)}_{\varepsilon} - \mathbf{s}_{\varepsilon} ) -
( \mathbf{s}_{\varepsilon} - \widetilde{\mathbf{u}}_0 ).
$$
Hence,
$$
\| \mathcal{A}_{\varepsilon} {\mathbf{v}}^{(1)}_{\varepsilon}
- \mathcal{A}^0 \widetilde{\mathbf{u}}_0 \|_{H^{-1} (\mathbb{R}^d)}  \leq
\| \mathcal{A}_{\varepsilon} ( {\mathbf{v}}^{(1)}_{\varepsilon} - \mathbf{s}_{\varepsilon} )
\|_{H^{-1} (\mathbb{R}^d)} + \| \mathbf{s}_{\varepsilon} - \widetilde{\mathbf{u}}_0
\|_{H^{-1} (\mathbb{R}^d)}.
\eqno(6.19)
$$
Next, taking (1.3) into account, we obtain
$$
\begin{aligned}
&\| \mathcal{A}_{\varepsilon} ({\mathbf{v}}^{(1)}_{\varepsilon} - {\mathbf{s}}_\varepsilon)
\|_{H^{-1} (\mathbb{R}^d)} = \sup_{0\ne \boldsymbol{\eta} \in H^1(\mathbb{R}^d;\mathbb{C}^n)}
\frac{\left|\left(g^\varepsilon b(\mathbf{D}) ({\mathbf{v}}^{(1)}_{\varepsilon} - {\mathbf{s}}_\varepsilon),
b(\mathbf{D}) \boldsymbol{\eta}\right)_{L_2(\mathbb{R}^d)}\right|}
{\|\boldsymbol{\eta}\|_{H^1(\mathbb{R}^d)}}
\\
&\leq \alpha_1 \|g\|_{L_\infty} \| {\mathbf{v}}^{(1)}_{\varepsilon}
- {\mathbf{s}}_\varepsilon\|_{H^1(\mathbb{R}^d)}.
\end{aligned}
$$
Combining this with (6.17)--(6.19), we see that
$$
\| \mathcal{A}_{\varepsilon} {\mathbf{v}}^{(1)}_{\varepsilon} -
\mathcal{A}^0 \widetilde{\mathbf{u}}_0 \|_{H^{-1} (\mathbb{R}^d)} \leq (C_1 + C_3 \alpha_1 \|g\|_{L_\infty})\varepsilon
\| \widetilde{\mathbf{F}}\|_{L_2(\mathbb{R}^d)}, \ \ 0 < \varepsilon \leq 1.
\eqno(6.20)
$$
Now, (6.15) and (6.20) imply (6.13) with the constant
$$
\widetilde{C}_4 =(C_1 + C_3 \alpha_1 \|g\|_{L_\infty})
\textrm{max} \left\{ \alpha_1 \| g \|_{L_\infty}, 1 \right\}.
$$

Returning to the case of a bounded domain, note that if
$\mathbf{f} \in H^{-1} (\mathcal{O}; \mathbb{C}^{n})$ and
$\widetilde{\mathbf{f}} \in H^{-1} (\mathbb{R}^d; \mathbb{C}^{n})$
are such that $\widetilde{\mathbf{f}} \vert_{\mathcal{O}} = \mathbf{f}$, then
$$
\begin{aligned}
&\| \mathbf{f} \|_{H^{-1} (\mathcal{O})} =
\sup_{0\ne \boldsymbol{\varphi} \in C_0^{\infty} (\mathcal{O})}
\frac{\left| \int_{\mathcal{O}} \left\langle \mathbf{f}, \boldsymbol{\varphi} \right\rangle
\, d \mathbf{x} \right|}{\| \boldsymbol{\varphi} \|_{H^1 (\mathcal{O})}}
= \sup_{0\ne \boldsymbol{\varphi} \in C_0^{\infty} (\mathcal{O})}
\frac{\left| \int_{\mathbb{R}^d} \langle \widetilde{\mathbf{f}}, \boldsymbol{\varphi}
\rangle \, d \mathbf{x} \right|}{\| \boldsymbol{\varphi} \|_{H^1 (\mathbb{R}^d)}}
\\
 &\leq
 \sup_{0\ne \boldsymbol{\varphi} \in C_0^{\infty} (\mathbb{R}^d)}
 \frac{\left| \int_{\mathbb{R}^d} \langle \widetilde{\mathbf{f}}, \boldsymbol{\varphi} \rangle
 \, d \mathbf{x} \right|}{\| \boldsymbol{\varphi} \|_{H^1 (\mathbb{R}^d)}}
 = \| \widetilde{\mathbf{f}} \|_{H^{-1} (\mathbb{R}^d)}.
\end{aligned}
$$
Hence,
$$
\| \mathcal{A}_{\varepsilon} \mathbf{\check{v}}_{\varepsilon}
- \mathcal{A}^0 \mathbf{u}_0 \|_{H^{-1} (\mathcal{O})} \leq
\| \mathcal{A}_{\varepsilon} {\mathbf{v}}^{(1)}_{\varepsilon}
- \mathcal{A}^0 \widetilde{\mathbf{u}}_0 \|_{H^{-1} (\mathbb{R}^d)}.
$$
Together with (6.13) and (6.10) this yields
$$
\| \mathcal{A}_{\varepsilon} \mathbf{\check{v}}_{\varepsilon}
- \mathcal{A}^0 \mathbf{u}_0 \|_{H^{-1} (\mathcal{O})} \leq
\widetilde{C}_4 \varepsilon \| \widetilde{\mathbf{u}}_0 \|_{H^2 (\mathbb{R}^d)}
\leq \widetilde{C}_4 C_{\mathcal{O}} \varepsilon \| \mathbf{u}_0 \|_{H^2 (\mathcal{O})}.
$$
Thus, inequality (6.12) holds with $C_4 = \widetilde{C}_4 C_\mathcal{O}$. \ $\bullet$

\smallskip
\noindent\textbf{6.3.}
The fist order approximation $\mathbf{\check{v}}_{\varepsilon}$ of the solution
$\mathbf{u}_\varepsilon$ defined by (6.2) does not satisfy the Dirichlet condition on
$\partial \mathcal{O}$. We consider the "discrepancy"\ $\mathbf{\check{w}}_{\varepsilon}$
which is the generalized solution of the problem
$$
\mathcal{A}_{\varepsilon} \mathbf{\check{w}}_{\varepsilon} = 0\ \ \text{in}\ \mathcal{O},
\quad \mathbf{\check{w}}_{\varepsilon} |_{\partial \mathcal{O}} =
\mathbf{\check{v}}_{\varepsilon} |_{\partial \mathcal{O}} =
\varepsilon \Lambda^{\varepsilon} b (\mathbf{D}) \mathbf{u}_0 |_{\partial \mathcal{O}}.
\eqno(6.21)
$$
Here the equation is understood in the weak sense: the function $\mathbf{\check{w}}_{\varepsilon} \in
H^1(\mathcal{O};\mathbb{C}^{n})$ satisfies the identity
$$
\int_{\mathcal{O}} \left\langle g^{\varepsilon} (\mathbf{x}) b (\mathbf{D})
\mathbf{\check{w}}_{\varepsilon}, b(\mathbf{D}) \boldsymbol{\eta} \right\rangle \, d \mathbf{x} = 0,
\quad \forall \boldsymbol{\eta} \in H^1_0 (\mathcal{O}; \mathbb{C}^{n}).
$$
The boundary condition in (6.21) is understood in the sense of the trace theorem:
under the assumptions of Theorem 6.1
one has $\Lambda^{\varepsilon} b (\mathbf{D}) \mathbf{u}_0
\in H^1 (\mathcal{O}; \mathbb{C}^{n})$, whence
$\Lambda^{\varepsilon} b (\mathbf{D}) \mathbf{u}_0 |_{\partial \mathcal{O}}
\in H^{1/2} (\partial \mathcal{O}; \mathbb{C}^{n})$.

By (4.2) and (4.18),
$\mathcal{A}_\varepsilon (\mathbf{u}_\varepsilon - \mathbf{\check{v}}_{\varepsilon})=
\mathcal{A}^0 \mathbf{u}_0 - \mathcal{A}_\varepsilon \mathbf{\check{v}}_{\varepsilon}$.
Consequently, by (6.21), the function
$\mathbf{u}_\varepsilon - \mathbf{\check{v}}_{\varepsilon} + \mathbf{\check{w}}_{\varepsilon}$
is the solution of the following Dirichlet problem
$$
\mathcal{A}_{\varepsilon} (\mathbf{u}_{\varepsilon} -
\mathbf{\check{v}}_{\varepsilon} + \mathbf{\check{w}}_{\varepsilon}) =
\mathcal{A}^0 \mathbf{u}_0 - \mathcal{A}_{\varepsilon} \mathbf{\check{v}}_{\varepsilon}
\ \ \text{in}\ \mathcal{O}, \quad
(\mathbf{u}_{\varepsilon} - \mathbf{\check{v}}_{\varepsilon} +
\mathbf{\check{w}}_{\varepsilon}) |_{\partial \mathcal{O}} = 0.
$$
The right-hand side in the equation belongs to $H^{-1}(\mathcal{O};\mathbb{C}^n)$.
Then, applying Lemmas 4.1 and 6.4, for $0< \varepsilon \leq 1$ we obtain
$$
 \| \mathbf{u}_{\varepsilon} - \mathbf{\check{v}}_{\varepsilon} +
 \mathbf{\check{w}}_{\varepsilon} \|_{H^1 (\mathcal{O}; \mathbb{C}^{n})}  \leq
 \widehat{C} \|
\mathcal{A}^0 \mathbf{u}_0 - \mathcal{A}_{\varepsilon} \mathbf{\check{v}}_{\varepsilon}
\|_{H^{-1} (\mathcal{O}; \mathbb{C}^{n})}
\leq \widehat{C} C_4 \varepsilon \| \mathbf{u}_0 \|_{H^2(\mathcal{O}; \mathbb{C}^{n})}.
$$
Together with (4.19) this implies that
$$
 \| \mathbf{u}_{\varepsilon} - \mathbf{\check{v}}_{\varepsilon} +
 \mathbf{\check{w}}_{\varepsilon} \|_{H^1 (\mathcal{O}; \mathbb{C}^{n})}  \leq
\widehat{C} C_4 \widehat{c}\, \varepsilon \|\mathbf{F}\|_{L_2(\mathcal{O};\mathbb{C}^n)},
\quad 0 < \varepsilon \leq 1.
\eqno(6.22)
$$
Therefore, the proof of estimate (6.3) from Theorem 6.1 is reduced to estimating of
$\mathbf{\check{w}}_{\varepsilon}$ in $H^1(\mathcal{O};\mathbb{C}^n)$.

Assume that $0< \varepsilon \leq \varepsilon_1$,
where the number $\varepsilon_1 \in (0,1]$ is defined in Lemma 5.2.
Fix a smooth cut-off function $\theta_{\varepsilon} (\mathbf{x})$ in $\mathbb{R}^d$
supported in the $\varepsilon$--vicinity of the boundary $\partial \mathcal{O}$
and such that
$$
\begin{aligned}
&\theta_{\varepsilon} \in C_0^{\infty} (\mathbb{R}^d), \quad
\text{supp}\,\theta_\varepsilon \subset (\partial \mathcal{O})_\varepsilon,
\quad 0 \leq \theta_{\varepsilon} (\mathbf{x}) \leq 1,
\\
&\theta_{\varepsilon} (\mathbf{x}) |_{\partial \mathcal{O}} = 1, \quad
\varepsilon \left| \nabla \theta_{\varepsilon} (\mathbf{x}) \right| \leq \kappa = \textrm{const}.
\end{aligned}
\eqno(6.23)
$$
Consider the following function in $\mathbb{R}^d$:
$$
\check{\boldsymbol{\phi}}_{\varepsilon} (\mathbf{x}) =
\varepsilon \theta_{\varepsilon} (\mathbf{x}) \Lambda^{\varepsilon} (\mathbf{x}) b (\mathbf{D})
\widetilde{\mathbf{u}}_0 (\mathbf{x}).
\eqno(6.24)
$$
Then $\check{\boldsymbol{\phi}}_{\varepsilon} \in H^1(\mathbb{R}^d;\mathbb{C}^n)$
and $\check{\boldsymbol{\phi}}_{\varepsilon}
\vert_{\partial \mathcal{O}} = \varepsilon \Lambda^{\varepsilon}
b (\mathbf{D}) \mathbf{u}_0\vert_{\partial \mathcal{O}}.$
The problem (6.21) can be rewritten as:
$\mathcal{A}_{\varepsilon} \mathbf{\check{w}}_{\varepsilon} = 0$ in $\mathcal{O}$,
$\mathbf{\check{w}}_{\varepsilon} |_{\partial \mathcal{O}} =
\check{\boldsymbol{\phi}}_{\varepsilon} |_{\partial \mathcal{O}}$.
Applying Lemma 4.3, we obtain
$$
\| \mathbf{\check{w}}_{\varepsilon} \|_{H^1 (\mathcal{O}; \mathbb{C}^{n})}
\leq \gamma_0 \| \check{\boldsymbol{\phi}}_{\varepsilon} \|_{H^1 (\mathcal{O}; \mathbb{C}^{n})}.
\eqno(6.25)
$$
Thus, the proof of the required estimate
for the norm of $\mathbf{\check{w}}_{\varepsilon}$ in
$H^1(\mathcal{O};\mathbb{C}^n)$ is reduced to the next statement.

\smallskip\noindent\textbf{Lemma 6.5.} \textit{Suppose that the assumptions of Theorem} 6.1
\textit{are satisfied. Assume that} \hbox{$0< \varepsilon \leq \varepsilon_1$},
\textit{where the number} $\varepsilon_1 \in (0,1]$ \textit{is defined in Lemma} 5.2.
\textit{Let} $\check{\boldsymbol{\phi}}_{\varepsilon}$ \textit{be the function defined
in accordance with} (6.23), (6.24). \textit{Then we have}
$$
\| \check{\boldsymbol{\phi}}_{\varepsilon} \|_{H^1 (\mathcal{O}; \mathbb{C}^{n})} \leq
C_5 \varepsilon^{1/2} \| \mathbf{F} \|_{L_2 (\mathcal{O}; \mathbb{C}^{n})},
\quad 0< \varepsilon \leq \varepsilon_1.
\eqno(6.26)
$$
\textit{The constant} $C_5$ \textit{depends only on} $m$, $d$, $\alpha_0$, $\alpha_1$,
$\| g \|_{L_\infty}$, $\| g^{-1} \|_{L_\infty}$, \textit{the norm}
$\| \Lambda \|_{L_\infty}$, \textit{and the domain} $\mathcal{O}$.

\smallskip\noindent\textbf{Proof.}
The norm of $\check{\boldsymbol{\phi}}_{\varepsilon}$ in $L_2(\mathcal{O};\mathbb{C}^n)$
is estimated with the help of Condition 1.9 and relations (4.17), (6.8), and (6.23):
$$
\begin{aligned}
&\| \check{\boldsymbol{\phi}}_{\varepsilon} \|_{L_2 (\mathcal{O})} \leq
\varepsilon  \| \Lambda^{\varepsilon} b (\mathbf{D}) \mathbf{u}_0 \|_{L_2 (\mathcal{O})}
\leq \varepsilon \| \Lambda \|_{L_\infty} \|b(\mathbf{D})\mathbf{u}_0 \|_{L_2(\mathcal{O})}
\\
&\leq \varepsilon \alpha_1^{1/2} \| \Lambda \|_{L_\infty}\| \mathbf{u}_0 \|_{H^1 (\mathcal{O})}
\leq \varepsilon \widehat{C} \alpha_1^{1/2} \| \Lambda \|_{L_\infty}\| \mathbf{F} \|_{L_2(\mathcal{O})}.
\end{aligned}
\eqno(6.27)
$$

Consider the derivatives
$$
\begin{aligned}
\frac{\partial \check{\boldsymbol{\phi}}_{\varepsilon}}{\partial x_j} =
\varepsilon \frac{\partial \theta_{\varepsilon}}{\partial x_j}
\Lambda^{\varepsilon} b (\mathbf{D}) \widetilde{\mathbf{u}}_0
+ \theta_{\varepsilon}
\left( \frac{\partial \Lambda}{\partial x_j}\right)^\varepsilon
b (\mathbf{D}) \widetilde{\mathbf{u}}_0
+ \varepsilon \theta_{\varepsilon} \Lambda^{\varepsilon}
(b (\mathbf{D}) \partial_j \widetilde{\mathbf{u}}_0),
\\
j=1,\dots,d.
\end{aligned}
$$
Then
$$
\begin{aligned}
\| \mathbf{D} \check{\boldsymbol{\phi}}_{\varepsilon}\|^2_{L_2(\mathcal{O})} &\leq
3 \varepsilon^2 \int_\mathcal{O} |\nabla \theta_\varepsilon |^2
|\Lambda^\varepsilon b(\mathbf{D})\mathbf{u}_0|^2\,d\mathbf{x} +
3 \int_\mathcal{O} |(\mathbf{D} \Lambda)^\varepsilon|^2
|\theta_\varepsilon b(\mathbf{D})\widetilde{\mathbf{u}}_0|^2\,d\mathbf{x}
\\
&+3 \varepsilon^2 \sum_{j=1}^d \int_\mathcal{O} |\theta_\varepsilon|^2
|\Lambda^\varepsilon b(\mathbf{D}) D_j \mathbf{u}_0|^2\,d\mathbf{x}.
\end{aligned}
\eqno(6.28)
$$
Denote the terms in the right-hand side of (6.28) by
${\mathcal I}_1$, ${\mathcal I}_2$, and ${\mathcal I}_3$, respectively.

It is easy to estimate ${\mathcal I}_3$.
By (6.23), Condition 1.9, and (1.2), (2.6), we obtain
$$
\left\| \theta_{\varepsilon} \Lambda^{\varepsilon}
b (\mathbf{D}) D_j \mathbf{u}_0 \right\|^2_{L_2(\mathcal{O})} \leq
\|\Lambda \|^2_{L_\infty} \alpha_1 d \sum_{l=1}^d \|D_l D_j \mathbf{u}_0\|^2_{L_2(\mathcal{O})}.
$$
Together with (4.19) this yields
$$
{\mathcal I}_3 \leq 3 \varepsilon^2 \|\Lambda \|^2_{L_\infty} \alpha_1 d
\|\mathbf{u}_0\|^2_{H^2(\mathcal{O})}
\leq \gamma_3 \varepsilon^2 \|\mathbf{F} \|^2_{L_2(\mathcal{O})},
\eqno(6.29)
$$
where $\gamma_3 = 3 \widehat{c}^{\,2} \alpha_1 d \| \Lambda \|^2_{L_\infty}$.

In order to estimate the first term in the right-hand side of (6.28), we apply (6.23),
Condition 1.9 and Lemma 5.1. We have
$$
\begin{aligned}
&{\mathcal I}_1 \leq 3 \kappa^2 \|\Lambda\|^2_{L_\infty}
\int_{B_\varepsilon} | b (\mathbf{D}) \mathbf{u}_0 |^2 \, d \mathbf{x}
\\
&\leq 3 \kappa^2 \| \Lambda \|^2_{L_\infty} \beta \varepsilon
\| b (\mathbf{D})\mathbf{u}_0 \|_{H^1 (\mathcal{O})} \| b(\mathbf{D}) \mathbf{u}_0 \|_{L_2(\mathcal{O})}.
\end{aligned}
$$
Using (4.17), (4.19), (6.8), and the estimate
$$
\| b(\mathbf{D}) \mathbf{u}_0 \|_{H^1(\mathcal{O})} =
\bigl\|\sum_{l=1}^d b_l D_l \mathbf{u}_0 \bigr\|_{H^1(\mathcal{O})}
\leq \alpha_1^{1/2} d^{1/2} \|\mathbf{u}_0\|_{H^2(\mathcal{O})},
$$
we arrive at the inequality
$$
{\mathcal I}_1 \leq 3 \varepsilon \kappa^2 \| \Lambda \|^2_{L_\infty} \beta
\alpha_1 d^{1/2}\|\mathbf{u}_0\|_{H^1(\mathcal{O})} \|\mathbf{u}_0\|_{H^2(\mathcal{O})}
\leq \gamma_1 \varepsilon \|\mathbf{F}\|^2_{L_2(\mathcal{O})},
\eqno(6.30)
$$
where $\gamma_1 = 3 \widehat{c} \widehat{C} \kappa^2
\| \Lambda \|^2_{L_\infty} \beta \alpha_1 d^{1/2}$.

It remains to consider the second term in the right-hand side of (6.28).
By Corollary 2.4,
$$
\begin{aligned}
&{\mathcal I}_2 \leq
3 \int_{\mathbb{R}^d} |(\mathbf{D} \Lambda)^\varepsilon|^2 |\theta_\varepsilon b(\mathbf{D})\widetilde{\mathbf{u}}_0|^2\,d\mathbf{x}
\\
&\leq 3\beta_1 \left\| \theta_{\varepsilon}
b (\mathbf{D}) \widetilde{\mathbf{u}}_0 \right\|^2_{L_2(\mathbb{R}^d)}
+ 3\beta_2 \|\Lambda\|_{L_\infty}^2 \varepsilon^2
\int_{\mathbb{R}^d} \sum_{l=1}^d \left| \frac{\partial}{\partial x_l}
(\theta_{\varepsilon} b (\mathbf{D}) \widetilde{\mathbf{u}}_0 ) \right|^2
\, d \mathbf{x}.
\end{aligned}
\eqno(6.31)
$$
Since
$$
\frac{\partial}{\partial x_l}
(\theta_{\varepsilon} b (\mathbf{D}) \widetilde{\mathbf{u}}_0)
=
\frac{\partial \theta_{\varepsilon}}{\partial x_l}
b (\mathbf{D}) \widetilde{\mathbf{u}}_0
+\theta_{\varepsilon} \frac{\partial}{\partial x_l}
(b (\mathbf{D}) \widetilde{\mathbf{u}}_0),
$$
then, by (6.23) and (6.31), we have
$$
\begin{aligned}
{\mathcal I}_2 &\leq 3 \left( \beta_1 + 2 \beta_2 \|\Lambda\|^2_{L_\infty} \kappa^2 \right)
 \int_{(\partial \mathcal{O})_{\varepsilon}} \left| b (\mathbf{D})
 \widetilde{\mathbf{u}}_0 \right|^2 \, d \mathbf{x}
\\
&+6 \beta_2 \|\Lambda\|^2_{L_\infty} \varepsilon^2
\int_{\mathbb{R}^d} \sum_{l=1}^d \left| b(\mathbf{D}) D_l \widetilde{\mathbf{u}}_0 \right|^2\,d\mathbf{x}.
\end{aligned}
$$
Combining this with Lemma 5.2 and condition (1.3), we obtain
$$
\begin{aligned}
{\mathcal I}_2 &\leq 3 \left( \beta_1 + 2 \beta_2 \|\Lambda\|^2_{L_\infty} \kappa^2 \right)
\beta^0 \varepsilon \|b(\mathbf{D}) \widetilde{\mathbf{u}}_0\|_{H^1(\mathbb{R}^d)}
\|b(\mathbf{D}) \widetilde{\mathbf{u}}_0\|_{L_2(\mathbb{R}^d)}
\\
&+6 \beta_2 \|\Lambda\|^2_{L_\infty} \varepsilon^2 \alpha_1
\| \widetilde{\mathbf{u}}_0\|^2_{H^2(\mathbb{R}^d)}
\\
&\leq
3 \varepsilon \left( \beta_1 + 2 \beta_2 \|\Lambda\|^2_{L_\infty} \kappa^2 \right) \beta^0 \alpha_1
\| \widetilde{\mathbf{u}}_0\|_{H^2(\mathbb{R}^d)}\| \widetilde{\mathbf{u}}_0\|_{H^1(\mathbb{R}^d)}
\\
&+ 6 \varepsilon^2 \beta_2 \|\Lambda\|^2_{L_\infty} \alpha_1
\| \widetilde{\mathbf{u}}_0\|^2_{H^2(\mathbb{R}^d)}.
\end{aligned}
\eqno(6.32)
$$
Taking (4.19) and (6.10) into account, from (6.32) we deduce that
$$
{\mathcal I}_2 \leq \gamma_2 \varepsilon \|\mathbf{F}\|^2_{L_2(\mathcal{O})},
\eqno(6.33)
$$
where $\gamma_2 = 3 (\widehat{c} C_\mathcal{O})^2 \left( \left( \beta_1 + 2 \beta_2
\|\Lambda\|^2_{L_\infty} \kappa^2 \right) \beta^0 \alpha_1
+ 2 \beta_2 \alpha_1 \|\Lambda\|^2_{L_\infty} \right)$.

Now, relations (6.28)--(6.30) and (6.33) imply that
$$
\left\| \mathbf{D} \check{\boldsymbol{\phi}}_{\varepsilon} \right\|^2_{L_2(\mathcal{O})}
\leq {\mathcal I}_1 + {\mathcal I}_2 + {\mathcal I}_3 \leq
\varepsilon(\gamma_1+ \gamma_2 + \gamma_3)\|\mathbf{F}\|^2_{L_2(\mathcal{O})},
\ \ 0< \varepsilon \leq \varepsilon_1.
\eqno(6.34)
$$
Finally, (6.27) and (6.34) yield estimate (6.26) with
$$
C_5=\left( \widehat{C}^2 \alpha_1 \|\Lambda \|^2_{L_\infty}
+\gamma_1+ \gamma_2 + \gamma_3 \right)^{1/2}. \ \ \bullet
$$

\smallskip

Now it is easy to complete the \textbf{proof of Theorem 6.1}. From (6.22),
(6.25), and (6.26) it follows that
$$
 \| \mathbf{u}_{\varepsilon} - \mathbf{\check{v}}_{\varepsilon}
 \|_{H^1 (\mathcal{O}; \mathbb{C}^{n})} \leq  \widehat{C} C_4 \widehat{c} \varepsilon
\|\mathbf{F}\|_{L_2(\mathcal{O};\mathbb{C}^n)} + \gamma_0 C_5 \varepsilon^{1/2} \|\mathbf{F}\|_{L_2(\mathcal{O};\mathbb{C}^n)},
\quad 0< \varepsilon \leq \varepsilon_1,
$$
which implies (6.3) with
$C_0 = \widehat{C} C_4 \widehat{c} + \gamma_0 C_5$.

It remains to check (6.4). From (6.3), (1.2), and (2.6) it follows that
$$
\begin{aligned}
&\| \mathbf{p}_\varepsilon - g^\varepsilon b(\mathbf{D})\mathbf{u}_0 -
g^\varepsilon b(\mathbf{D})(\varepsilon \Lambda^\varepsilon
b(\mathbf{D})\mathbf{u}_0) \|_{L_2(\mathcal{O})}
\\
&\leq \|g\|_{L_\infty} \alpha_1^{1/2} d^{1/2}C_0 \varepsilon^{1/2} \|\mathbf{F}\|_{L_2(\mathcal{O})},
\ \ 0< \varepsilon \leq \varepsilon_1.
\end{aligned}
\eqno(6.35)
$$
From (1.2) and the definition of the matrix $\widetilde{g}$ it is seen that
$$
g^\varepsilon b(\mathbf{D})\mathbf{u}_0 + g^\varepsilon b(\mathbf{D})(\varepsilon \Lambda^\varepsilon
b(\mathbf{D})\mathbf{u}_0) = \widetilde{g}^\varepsilon b(\mathbf{D}) \mathbf{u}_0
+ \varepsilon g^\varepsilon \sum_{l=1}^d b_l \Lambda^\varepsilon b(\mathbf{D}) D_l \mathbf{u}_0.
\eqno(6.36)
$$
Applying Condition 1.9 and relations (1.2), (2.6), and (4.19), we obtain
$$
\bigl\|\varepsilon g^\varepsilon \sum_{l=1}^d b_l \Lambda^\varepsilon b(\mathbf{D}) D_l \mathbf{u}_0
\bigr\|_{L_2(\mathcal{O})} \leq \varepsilon \|g\|_{L_\infty} \|\Lambda\|_{L_\infty}
\alpha_1 d \widehat{c} \|\mathbf{F}\|_{L_2(\mathcal{O})}.
\eqno(6.37)
$$
Now, relations (6.35)--(6.37) imply (6.4) with the constant
$C_0'= \|g\|_{L_\infty} \alpha_1^{1/2} d^{1/2}C_0 + \|g\|_{L_\infty} \|\Lambda\|_{L_\infty}
\alpha_1 d \widehat{c}$. \ $\bullet$

\section*{\S 7. Results in the general case}

\noindent\textbf{7.1.} Now we refuse the assumption that $\Lambda(\mathbf{x})$ is bounded.
Then we need to include a smoothing operator in the corrector.

Let $P_\mathcal{O}$ be the extension operator (6.9), and let $S_\varepsilon$
be the operator smoothing in Steklov's sense defined by (3.1).
By $R_\mathcal{O}$ we denote the operator of restriction of functions in $\mathbb{R}^d$
to the domain $\mathcal{O}$. We put
$$
K_D(\varepsilon) = R_\mathcal{O} [\Lambda^\varepsilon] S_\varepsilon b(\mathbf{D})
P_\mathcal{O} (\mathcal{A}_D^0)^{-1}.
\eqno(7.1)
$$
The operator $b(\mathbf{D}) P_\mathcal{O} (\mathcal{A}_D^0)^{-1}$ is a continuous mapping of
$L_2(\mathcal{O};\mathbb{C}^n)$ into $H^1(\mathbb{R}^d;\mathbb{C}^m)$.
As has been mentioned in Subsection 3.2, the operator $[\Lambda^\varepsilon] S_\varepsilon$
is continuous from $H^1(\mathbb{R}^d;\mathbb{C}^m)$ to $H^1(\mathbb{R}^d;\mathbb{C}^n)$.
Consequently, the operator (7.1) is continuous from
$L_2(\mathcal{O};\mathbb{C}^n)$ to $H^1(\mathcal{O};\mathbb{C}^n)$.

Let $\mathbf{u}_\varepsilon$ be the solution of the problem (4.2),
and let $\mathbf{u}_0$ be the solution of the problem (4.18).
As above, we denote $\widetilde{\mathbf{u}}_0 = P_\mathcal{O} \mathbf{u}_0$.
We put
$$
{\mathbf{v}}^{(2)}_\varepsilon(\mathbf{x}) = \widetilde{\mathbf{u}}_0(\mathbf{x})
+ \varepsilon \Lambda^\varepsilon(\mathbf{x}) (S_\varepsilon b(\mathbf{D})
\widetilde{\mathbf{u}}_0)(\mathbf{x}),
$$
and $\mathbf{v}_\varepsilon := {\mathbf{v}}^{(2)}_\varepsilon \vert_{\mathcal{O}}$. Then
$$
\mathbf{v}_\varepsilon = (\mathcal{A}_D^0)^{-1} \mathbf{F} + \varepsilon K_D(\varepsilon) \mathbf{F}.
\eqno(7.2)
$$

The following theorem is \textit{our main result} in the general case.

\smallskip\noindent\textbf{Theorem 7.1.} \textit{Suppose that}
$\mathcal{O}\subset \mathbb{R}^d$ \textit{is a bounded domain of class} $C^2$.
\textit{Let} $g(\mathbf{x})$ \textit{and} $b(\mathbf{D})$ \textit{satisfy the assumptions
of Subsection} 1.2. \textit{Let} $\mathbf{u}_{\varepsilon}$ \textit{be the solution of the problem} (4.2),
\textit{and let} $\mathbf{u}_0$ \textit{be the solution of the problem задачи} (4.18)
\textit{with} $\mathbf{F} \in L_2(\mathcal{O};\mathbb{C}^n)$.
\textit{Let} $\mathbf{v}_{\varepsilon}$ \textit{be the function defined by} (7.1), (7.2).
\textit{Then there exists a number} $\varepsilon_2 \in (0,1]$
\textit{depending on the domain} $\mathcal{O}$ \textit{and the lattice} $\Gamma$
\textit{such that we have}
$$
\| \mathbf{u}_{\varepsilon} - \mathbf{v}_{\varepsilon}
\|_{H^1 (\mathcal{O}; \mathbb{C}^{n})} \leq C \varepsilon^{1/2}
\| \mathbf{F} \|_{L_2 (\mathcal{O}; \mathbb{C}^{n})},\quad 0< \varepsilon \leq \varepsilon_2,
\eqno(7.3)
$$
\textit{or, in operator terms,}
$$
\| \mathcal{A}_{D,\varepsilon}^{-1} - (\mathcal{A}^0_D)^{-1} -
\varepsilon K_D (\varepsilon) \|_{L_2 (\mathcal{O}; \mathbb{C}^{n}) \to
H^1 (\mathcal{O}; \mathbb{C}^{n})} \leq C {\varepsilon}^{1/2}.
$$
\textit{The flux} $\mathbf{p}_\varepsilon:= g^\varepsilon b(\mathbf{D})\mathbf{u}_\varepsilon$
\textit{admits the following approximation}
$$
\| \mathbf{p}_\varepsilon - \widetilde{g}^\varepsilon S_\varepsilon b(\mathbf{D}) \widetilde{\mathbf{u}}_0
\|_{L_2(\mathcal{O};\mathbb{C}^m)} \leq C' \varepsilon^{1/2}
\| \mathbf{F} \|_{L_2 (\mathcal{O}; \mathbb{C}^{n})},\ \ 0< \varepsilon \leq \varepsilon_2,
\eqno(7.4)
$$
\textit{where} $\widetilde{g}(\mathbf{x}):= g(\mathbf{x})(b(\mathbf{D})\Lambda(\mathbf{x}) + \mathbf{1}_m)$.
\textit{The constants} $C$, $C'$ \textit{depend only on}
$m$, $d$, $\alpha_0$, $\alpha_1$, $\| g \|_{L_\infty}$,
$\| g^{-1} \|_{L_\infty}$, \textit{the parameters of the lattice} $\Gamma$,
\textit{and the domain} $\mathcal{O}$.

Roughening the result of Theorem 7.1, we arrive at the following corollary.

\smallskip\noindent\textbf{Corollary 7.2.}
\textit{Under the assumptions of Theorem} 7.1 \textit{for} $0< \varepsilon \leq \varepsilon_2$
\textit{we have}
$$
\| \mathbf{u}_{\varepsilon} - \mathbf{u}_0
\|_{L_2 (\mathcal{O}; \mathbb{C}^{n})} \leq \widetilde{C} \varepsilon^{1/2}
\| \mathbf{F} \|_{L_2 (\mathcal{O}; \mathbb{C}^{n})},
\eqno(7.5)
$$
\textit{or, in operator terms,}
$$
\| \mathcal{A}_{D,\varepsilon}^{-1} - (\mathcal{A}^0_D)^{-1}
\|_{L_2 (\mathcal{O}; \mathbb{C}^{n}) \to
L_2(\mathcal{O}; \mathbb{C}^{n})} \leq \widetilde{C} {\varepsilon}^{1/2}.
$$
\textit{The constant} $\widetilde{C}$ \textit{is given by}
$$
\widetilde{C} = C + C_\mathcal{O} \widehat{c} m^{1/2} (2r_0)^{-1}
\alpha_0^{-1/2} \alpha_1^{1/2} \|g\|^{1/2}_{L_\infty} \|g^{-1}\|^{1/2}_{L_\infty},
$$
\textit{where} $\widehat{c}$ \textit{is the constant from} (4.19),
$C_\mathcal{O}$ \textit{is the norm of the extension operator}
$P_\mathcal{O}$, \textit{and} $r_0$ \textit{is the radius of the ball inscribed in}
$\text{clos} \, \widetilde{\Omega}$.

\smallskip\noindent\textbf{Proof.} From (7.2) and (7.3) it follows that
$$
\| \mathbf{u}_\varepsilon - \mathbf{u}_0 \|_{L_2(\mathcal{O})}
\le C \varepsilon^{1/2} \|\mathbf{F}\|_{L_2(\mathcal{O})}+
\varepsilon \| \Lambda^\varepsilon S_\varepsilon b(\mathbf{D})
\widetilde{\mathbf{u}}_0 \|_{L_2(\mathcal{O})}.
\eqno(7.6)
$$
By (3.9) and (1.3),
$$
\begin{aligned}
\| \Lambda^\varepsilon S_\varepsilon b(\mathbf{D}) \widetilde{\mathbf{u}}_0 \|_{L_2(\mathcal{O})} \leq
\| \Lambda^\varepsilon S_\varepsilon b(\mathbf{D}) \widetilde{\mathbf{u}}_0 \|_{L_2(\mathbb{R}^d)}
\\
\leq M \| b(\mathbf{D}) \widetilde{\mathbf{u}}_0 \|_{L_2(\mathbb{R}^d)}
\leq M \alpha_1^{1/2} \|\widetilde{\mathbf{u}}_0 \|_{H^1(\mathbb{R}^d)}.
\end{aligned}
\eqno(7.7)
$$
Taking (4.19) and (6.10) into account, we obtain
$$
\| \widetilde{\mathbf{u}}_0 \|_{H^1(\mathbb{R}^d)} \leq \| \widetilde{\mathbf{u}}_0 \|_{H^2(\mathbb{R}^d)}
\leq C_\mathcal{O} \|\mathbf{u}_0\|_{H^2(\mathcal{O})} \leq C_\mathcal{O} \widehat{c}
\|\mathbf{F} \|_{L_2(\mathcal{O})}.
\eqno(7.8)
$$
Now, from (7.6)--(7.8) it follows that
$$
\| \mathbf{u}_\varepsilon - \mathbf{u}_0 \|_{L_2(\mathcal{O})} \leq
C \varepsilon^{1/2} \|\mathbf{F}\|_{L_2(\mathcal{O})}+
\varepsilon M \alpha_1^{1/2}C_\mathcal{O} \widehat{c} \|\mathbf{F} \|_{L_2(\mathcal{O})}.
$$
Recalling the expression for $M$ (see (3.8)), we arrive at (7.5). $\bullet$

\smallskip\noindent\textbf{7.2.} Let us start the proof of Theorem 7.1.
The following statement is similar to Lemma 6.4.

\smallskip\noindent\textbf{Lemma 7.3.} \textit{Let} $\mathbf{u}_0$
\textit{be the solution of the problem} (4.18), \textit{and let} $\mathbf{v}_\varepsilon$
\textit{be the function defined by} (7.1), (7.2). \textit{Then for} $0< \varepsilon \leq 1$
\textit{we have}
$$
\| \mathcal{A}_{\varepsilon} \mathbf{v}_{\varepsilon} -
\mathcal{A}^0 \mathbf{u}_0 \|_{H^{-1} (\mathcal{O}; \mathbb{C}^{n})}
\leq  C_6 \varepsilon \| \mathbf{u}_0 \|_{H^2 (\mathcal{O}; \mathbb{C}^{n})}.
$$
\textit{Here the constant} $C_6$ \textit{is given by}
$$
C_6= C_\mathcal{O} (C_1 + \widetilde{C}_2 \alpha_1 \|g\|_{L_\infty}) \max \{\alpha_1 \|g\|_{L_\infty},1\}
$$
\textit{and depends only on} $m$, $d$, $\alpha_0$, $\alpha_1$, $\| g \|_{L_\infty}$,
$\| g^{-1} \|_{L_\infty}$, \textit{the parameters of the lattice} $\Gamma$,
\textit{and the domain} $\mathcal{O}$.

\smallskip\noindent\textbf{Proof.} Lemma 7.3 can be proved by analogy with the proof of Lemma 6.4.
The only difference is that one should apply Theorem 3.3 instead of Theorem 1.10. $\bullet$

\smallskip
Next, by analogy with the proof of Theorem 6.1, we consider the "discrepancy"\
$\mathbf{w}_\varepsilon \in H^1(\mathcal{O};\mathbb{C}^n)$ which is the generalized solution
of the problem
$$
\mathcal{A}_{\varepsilon} \mathbf{w}_{\varepsilon} = 0\ \ \text{in}\ \mathcal{O},
\quad \mathbf{w}_{\varepsilon} |_{\partial \mathcal{O}} =
\mathbf{v}_{\varepsilon} |_{\partial \mathcal{O}} =
\varepsilon \Lambda^{\varepsilon} (S_{\varepsilon} b (\mathbf{D})
\widetilde{\mathbf{u}}_0)\vert_{\partial \mathcal{O}}.
\eqno(7.9)
$$
The equation in (7.9) is understood in the weak sense,
and the boundary condition in the sense of the trace theorem.
It should be taken into account that $\Lambda^\varepsilon (S_{\varepsilon} b (\mathbf{D})
\widetilde{\mathbf{u}}_0) \in H^1(\mathcal{O};\mathbb{C}^n)$.

By (4.2), (4.18), and (7.9), the function $\mathbf{u}_\varepsilon - \mathbf{v}_\varepsilon +
\mathbf{w}_\varepsilon$ is the solution of the following problem
$$
\mathcal{A}_\varepsilon (\mathbf{u}_\varepsilon - \mathbf{v}_\varepsilon + \mathbf{w}_\varepsilon)
= \mathcal{A}^0 \mathbf{u}_0 - \mathcal{A}_\varepsilon \mathbf{v}_\varepsilon \ \ \text{in}\ \mathcal{O},
\quad (\mathbf{u}_\varepsilon - \mathbf{v}_\varepsilon + \mathbf{w}_\varepsilon)\vert_{\partial \mathcal{O}} =0.
$$

Applying Lemmas 4.1 and 7.3, for $0< \varepsilon \leq 1$ we obtain
$$
\| \mathbf{u}_\varepsilon - \mathbf{v}_\varepsilon + \mathbf{w}_\varepsilon
\|_{H^1(\mathcal{O};\mathbb{C}^n)} \leq \widehat{C}
\| \mathcal{A}^0 \mathbf{u}_0 - \mathcal{A}_{\varepsilon} \mathbf{v}_{\varepsilon}
\|_{H^{-1} (\mathcal{O}; \mathbb{C}^{n})} \leq \widehat{C} C_6 \varepsilon \| \mathbf{u}_0 \|_{H^2(\mathcal{O};\mathbb{C}^n)}.
$$
Together with (4.19) this implies that
$$
\| \mathbf{u}_\varepsilon - \mathbf{v}_\varepsilon + \mathbf{w}_\varepsilon
\|_{H^1(\mathcal{O};\mathbb{C}^n)} \leq \widehat{C} C_6 \widehat{c} \varepsilon
\| \mathbf{F} \|_{L_2(\mathcal{O};\mathbb{C}^n)}, \ \ 0< \varepsilon \leq 1.
\eqno(7.10)
$$

\smallskip\noindent\textbf{7.3.} By (7.10), the proof of estimate (7.3)
is reduced to estimating of the $H^1$-norm of $\mathbf{w}_\varepsilon$.
As in Subsection 6.3, we fix a cut-off function $\theta_\varepsilon(\mathbf{x})$
satisfying conditions (6.23).
We assume that $0< \varepsilon \leq \varepsilon_2$, where the number
$\varepsilon_2 \in (0,1]$ is defined in Lemma 5.3.
Consider the following function in $\mathbb{R}^d$:
$$
{\boldsymbol{\phi}}_\varepsilon(\mathbf{x})= \varepsilon \theta_\varepsilon(\mathbf{x})
\Lambda^\varepsilon(\mathbf{x}) (S_\varepsilon b(\mathbf{D}) \widetilde{\mathbf{u}}_0)(\mathbf{x}).
\eqno(7.11)
$$
Similarly to (6.25), by Lemma 4.3, we have
$$
\| \mathbf{w}_\varepsilon \|_{H^1(\mathcal{O};\mathbb{C}^n)} \leq
\gamma_0 \| {\boldsymbol{\phi}}_\varepsilon\|_{H^1(\mathcal{O};\mathbb{C}^n)}.
\eqno(7.12)
$$
Thus, the problem is reduced to the proof of the following statement.

\smallskip\noindent\textbf{Lemma 7.4.} \textit{Suppose that the assumptions of Theorem} 7.1
\textit{are satisfied. Let} \hbox{$0< \varepsilon \leq \varepsilon_2$},
\textit{where the number} $\varepsilon_2 \in (0,1]$ \textit{is defined in Lemma} 5.3.
\textit{Let} ${\boldsymbol{\phi}}_{\varepsilon}$ \textit{be the function defined
in accordance with} (6.23), (7.11). \textit{Then we have}
$$
\| {\boldsymbol{\phi}}_{\varepsilon} \|_{H^1 (\mathcal{O}; \mathbb{C}^{n})} \leq
C_7 \varepsilon^{1/2} \| \mathbf{F} \|_{L_2 (\mathcal{O}; \mathbb{C}^{n})},
\quad 0< \varepsilon \leq \varepsilon_2.
\eqno(7.13)
$$
\textit{The constant} $C_7$ \textit{depends only on} $m$, $d$, $\alpha_0$, $\alpha_1$,
$\|g\|_{L_\infty}$, $\|g^{-1} \|_{L_\infty}$, \textit{the parameters of the lattice} $\Gamma$,
 \textit{and the domain} $\mathcal{O}$.

\smallskip\noindent\textbf{Proof.}
We start with the estimate for the norm of the function (7.11) in $L_2(\mathcal{O};\mathbb{C}^n)$.
From (1.3), (3.9), (6.23), and (7.8) it follows that
$$
\begin{aligned}
&\| {\boldsymbol{\phi}}_{\varepsilon} \|_{L_2 (\mathcal{O})} \leq
\varepsilon  \| \Lambda^{\varepsilon} (S_\varepsilon b (\mathbf{D})
\widetilde{\mathbf{u}}_0) \|_{L_2 (\mathbb{R}^d)}
\leq \varepsilon M \alpha_1^{1/2} \| \widetilde{\mathbf{u}}_0 \|_{H^1(\mathbb{R}^d)}
\\
&\leq \varepsilon M \alpha_1^{1/2} C_\mathcal{O} \widehat{c} \| \mathbf{F} \|_{L_2(\mathcal{O})}.
\end{aligned}
\eqno(7.14)
$$

Consider the derivatives
$$
\begin{aligned}
&\frac{\partial {\boldsymbol{\phi}}_{\varepsilon}}{\partial x_j} =
\varepsilon \frac{\partial \theta_{\varepsilon}}{\partial x_j}
\Lambda^{\varepsilon} (S_\varepsilon b (\mathbf{D}) \widetilde{\mathbf{u}}_0)
+ \theta_{\varepsilon}
\left( \frac{\partial \Lambda}{\partial x_j}\right)^\varepsilon
(S_\varepsilon b (\mathbf{D}) \widetilde{\mathbf{u}}_0)
\\
&+ \varepsilon \theta_{\varepsilon} \Lambda^{\varepsilon}
 (S_\varepsilon b (\mathbf{D}) \partial_j \widetilde{\mathbf{u}}_0),\ \ j=1,\dots,d.
\end{aligned}
$$
Then
$$
\begin{aligned}
&\| \mathbf{D} {\boldsymbol{\phi}}_{\varepsilon}\|^2_{L_2(\mathcal{O})} \leq
3 \varepsilon^2 \int_\mathcal{O} |\nabla \theta_\varepsilon |^2 |\Lambda^\varepsilon (S_\varepsilon b(\mathbf{D})\widetilde{\mathbf{u}}_0)|^2\,d\mathbf{x}
\\
&+
3 \int_\mathcal{O} |(\mathbf{D} \Lambda)^\varepsilon|^2 |\theta_\varepsilon (S_\varepsilon b(\mathbf{D})\widetilde{\mathbf{u}}_0)|^2\,d\mathbf{x}
+3 \varepsilon^2 \sum_{j=1}^d \int_\mathcal{O} |\theta_\varepsilon|^2 |\Lambda^\varepsilon (S_\varepsilon b(\mathbf{D}) D_j \widetilde{\mathbf{u}}_0)|^2\,d\mathbf{x}.
\end{aligned}
\eqno(7.15)
$$
The summands in the right-hand side of (7.15) are denoted by
${\mathcal J}_1$, ${\mathcal J}_2$, and ${\mathcal J}_3$, respectively.

It is easy to estimate ${\mathcal J}_3$. From
(1.3), (3.9), and (6.23) it follows that
$$
{\mathcal J}_3 \leq
3 \varepsilon^2 \sum_{j=1}^d \left\| \Lambda^{\varepsilon}
(S_\varepsilon b (\mathbf{D}) D_j \widetilde{\mathbf{u}}_0) \right\|^2_{L_2(\mathbb{R}^d)}
\leq 3 \varepsilon^2 M^2  \alpha_1 \| \widetilde{\mathbf{u}}_0\|^2_{H^2(\mathbb{R}^d)}.
$$
Combining this with (4.19) and (6.10), we obtain
$$
{\mathcal J}_3 \leq \widehat{\gamma}_3 \varepsilon^2 \|\mathbf{F} \|^2_{L_2(\mathcal{O})},
\eqno(7.16)
$$
where $\widehat{\gamma}_3 = 3 M^2 \alpha_1 (C_\mathcal{O} \widehat{c})^{2}$.

The fist term in the right-hand side of (7.15) is estimated with the help of
(6.23) and Lemma 5.3. For $0< \varepsilon \leq \varepsilon_2$ we have
$$
\begin{aligned}
&{\mathcal J}_1 \leq 3 \kappa^2 \int_{(\partial \mathcal{O})_\varepsilon} |\Lambda^\varepsilon
\left( S_\varepsilon b (\mathbf{D}) \widetilde{\mathbf{u}}_0\right)|^2\,d\mathbf{x}
\\
&\leq 3 \kappa^2 \beta_* \varepsilon |\Omega|^{-1}\|\Lambda\|^2_{L_2(\Omega)}
\| b (\mathbf{D}) \widetilde{\mathbf{u}}_0 \|_{H^1(\mathbb{R}^d)}
\| b (\mathbf{D}) \widetilde{\mathbf{u}}_0 \|_{L_2(\mathbb{R}^d)}.
\end{aligned}
$$
Combining this with (1.3), (4.19), (6.10), and estimate (2.14), we arrive at the inequality
$$
{\mathcal J}_1 \leq \widehat{\gamma}_1 \varepsilon \|\mathbf{F} \|^2_{L_2(\mathcal{O})},
\eqno(7.17)
$$
where $\widehat{\gamma}_1 = 3 \kappa^2 \beta_* (C_\mathcal{O} \widehat{c})^{2} m (2r_0)^{-2}
\alpha_0^{-1} \alpha_1 \|g\|_{L_\infty}\|g^{-1}\|_{L_\infty}$.

It remains to consider the second term in the right-hand side of (7.15). By (6.23),
$$
{\mathcal J}_2
\leq 3 \int_{(\partial \mathcal{O})_\varepsilon} |(\mathbf{D} \Lambda)^\varepsilon|^2
|S_\varepsilon b(\mathbf{D})\widetilde{\mathbf{u}}_0|^2\,d\mathbf{x}.
$$
By Lemma 5.3, for $0< \varepsilon \leq \varepsilon_2$ we have
$$
{\mathcal J}_2 \leq 3 \beta_* \varepsilon |\Omega|^{-1} \|\mathbf{D} \Lambda\|^2_{L_2(\Omega)}
\| b (\mathbf{D}) \widetilde{\mathbf{u}}_0 \|_{H^1(\mathbb{R}^d)}
\| b (\mathbf{D}) \widetilde{\mathbf{u}}_0 \|_{L_2(\mathbb{R}^d)}.
$$
Together with (1.3), (2.15), (4.19), and (6.10) this implies that
$$
{\mathcal J}_2 \leq \widehat{\gamma}_2 \varepsilon \| \mathbf{F} \|^2_{L_2(\mathcal{O})},
\ \ 0< \varepsilon \leq \varepsilon_2,
\eqno(7.18)
$$
where $\widehat{\gamma}_2 = 3 \beta_* (C_\mathcal{O} \widehat{c})^{2} m
\alpha_0^{-1} \alpha_1\|g\|_{L_\infty}\|g^{-1}\|_{L_\infty}$.

Finally, relations (7.15)--(7.18) yield
$$
\| \mathbf{D} {\boldsymbol{\phi}}_{\varepsilon}\|^2_{L_2(\mathcal{O})} \leq
{\mathcal J}_1 + {\mathcal J}_2 + {\mathcal J}_3 \leq
(\widehat{\gamma}_1 + \widehat{\gamma}_2 + \widehat{\gamma}_3) \varepsilon
\| \mathbf{F} \|^2_{L_2(\mathcal{O})},\ \ 0< \varepsilon \leq \varepsilon_2.
$$
Combining this with (7.14), we obtain (7.13) with
$$
C_7 =(M^2 \alpha_1 (C_\mathcal{O} \widehat{c})^{2} + \widehat{\gamma}_1
+ \widehat{\gamma}_2 + \widehat{\gamma}_3)^{1/2}.  \ \ \bullet
$$

\smallskip
Now, it is easy to complete the \textbf{proof of Theorem 7.1}.
From (7.10), (7.12), and (7.13) it follows that
$$
 \| \mathbf{u}_{\varepsilon} - \mathbf{v}_{\varepsilon}
 \|_{H^1 (\mathcal{O}; \mathbb{C}^{n})} \leq  \widehat{C} C_6 \widehat{c} \varepsilon
\|\mathbf{F}\|_{L_2(\mathcal{O};\mathbb{C}^n)} + \gamma_0 C_7 \varepsilon^{1/2}
\|\mathbf{F}\|_{L_2(\mathcal{O};\mathbb{C}^n)}, \quad 0< \varepsilon \leq \varepsilon_2.
$$
This implies (7.3) with $C = \widehat{C} C_6 \widehat{c} + \gamma_0 C_7$.

It remains to check (7.4). Taking (1.2) and (2.6) into account, from (7.3) we obtain
$$
\begin{aligned}
&\| \mathbf{p}_\varepsilon - g^\varepsilon b(\mathbf{D})\mathbf{u}_0 - g^\varepsilon b(\mathbf{D})(\varepsilon \Lambda^\varepsilon S_\varepsilon b(\mathbf{D})\widetilde{\mathbf{u}}_0) \|_{L_2(\mathcal{O})}
\\
&\leq \|g\|_{L_\infty} \alpha_1^{1/2} d^{1/2} C \varepsilon^{1/2} \|\mathbf{F}\|_{L_2(\mathcal{O})},
\ \ 0< \varepsilon \leq \varepsilon_2.
\end{aligned}
\eqno(7.19)
$$
By Proposition 3.1 and relations (1.3), (4.19), and (6.10), we conclude that
$$
\begin{aligned}
\| g^\varepsilon b(\mathbf{D})\mathbf{u}_0 - g^\varepsilon S_\varepsilon b(\mathbf{D})
\widetilde{\mathbf{u}}_0 \|_{L_2(\mathcal{O})} \leq
\|g\|_{L_\infty} \| b(\mathbf{D})\widetilde{\mathbf{u}}_0 -  S_\varepsilon b(\mathbf{D})
\widetilde{\mathbf{u}}_0 \|_{L_2(\mathbb{R}^d)}
\\
\leq \varepsilon r_1 \|g\|_{L_\infty} \alpha_1^{1/2} \|\widetilde{\mathbf{u}}_0\|_{H^2(\mathbb{R}^d)}
\leq \varepsilon r_1 \|g\|_{L_\infty} \alpha_1^{1/2} C_\mathcal{O}
\widehat{c} \|\mathbf{F}\|_{L_2(\mathcal{O})}.
\end{aligned}
\eqno(7.20)
$$
From (1.2) and the definition of the matrix $\widetilde{g}$ it is seen that
$$
\begin{aligned}
&g^\varepsilon S_\varepsilon b(\mathbf{D}) \widetilde{\mathbf{u}}_0 +
g^\varepsilon b(\mathbf{D})(\varepsilon \Lambda^\varepsilon S_\varepsilon b(\mathbf{D})
\widetilde{\mathbf{u}}_0)
\\
&= \widetilde{g}^\varepsilon S_\varepsilon b(\mathbf{D}) \widetilde{\mathbf{u}}_0
+ \varepsilon g^\varepsilon \sum_{l=1}^d b_l \Lambda^\varepsilon S_\varepsilon b(\mathbf{D})
D_l \widetilde{\mathbf{u}}_0.
\end{aligned}
\eqno(7.21)
$$
Taking (1.3), (2.6), (3.9), (4.19), and (6.10) into account, we obtain
$$
\begin{aligned}
\bigl\|\varepsilon g^\varepsilon \sum_{l=1}^d b_l \Lambda^\varepsilon S_\varepsilon
b(\mathbf{D}) D_l \widetilde{\mathbf{u}}_0 \bigr\|_{L_2(\mathcal{O})}
\leq \varepsilon \|g\|_{L_\infty} M \alpha_1 d^{1/2} \| \widetilde{\mathbf{u}}_0 \|_{H^2(\mathbb{R}^d)}
\\
\leq \varepsilon \|g\|_{L_\infty} M \alpha_1 d^{1/2} C_\mathcal{O} \widehat{c}
\|\mathbf{F}\|_{L_2(\mathcal{O})}.
\end{aligned}
\eqno(7.22)
$$
Now, relations (7.19)--(7.22) imply (7.4) with the constant
$C'= \|g\|_{L_\infty} \alpha_1^{1/2} d^{1/2} C
+ \|g\|_{L_\infty}  C_\mathcal{O} \widehat{c} (r_1 \alpha_1^{1/2} + M \alpha_1 d^{1/2}).$ \ $\bullet$


\begin{thebibliography}{13}

\bibitem[BaPa]{BaPa}Bakhvalov N.~S., Panasenko G.~P., \textit{Homogenization: averaging processes
in periodic media. Mathematical problems in mechanics of composite materials}, "Nauka", Moscow, 1984;
English transl., Math. Appl. (Soviet Ser.), vol. 36, Kluwer Acad. Publ. Group, Dordrecht, 1989.

\bibitem[BeLP]{BeLP} Bensoussan~A., Lions~J.-L., Papanicolaou~G.,
\emph{Asymptotic analysis for periodic structures}, Stud. Math. Appl., vol. 5,
North-Holland Publishing Co., Amsterdam-New York, 1978.

\bibitem[BSu1]{BSu1}Birman M.~Sh., Suslina T.~A., \textit{Threshold effects near the lower edge
of the spectrum for periodic differential operators of mathematical physics},
Systems, Approximation, Singular Integral Operators, and Related Topics (Bodeaux, 2000),
Oper. Theory Adv. Appl., vol. 129, Bikh\"auser, Basel, 2001, pp. 71--107.



\bibitem[BSu2]{BSu2}Birman M.~Sh., Suslina T.~A., \textit{Second order periodic
differential operators. Threshold properties and homogenization}, Algebra i Analiz
{\bf 15} (2003), no. 5, 1-108; English transl., St.~Petersburg Math. J. {\bf 15} (2004),
no. 5, 639--714.

\bibitem[BSu3]{BSu3}Birman M.~Sh., Suslina T.~A.,
\textit{Theshold approximations with corrector for the resolvent of a factorized selfadjoint
operator family}, Algebra i Analiz {\bf 17} (2005), no.~5, 69--90;
English transl., St.~Petersburg Math. J. {\bf 17} (2006), no.~5, 745--762.

\bibitem[BSu4]{BSu4}Birman M.~Sh., Suslina T.~A.,
\textit{Homogenization with corrector term for periodic elliptic differential
operators}, Algebra i Analiz {\bf 17} (2005), no.~6, 1--104;
English transl., St.~Petersburg Math. J. {\bf 17} (2006), no.~6, 897--973.

\bibitem[BSu5]{BSu5}Birman M.~Sh., Suslina T.~A.,
\textit{Homogenization with corrector term for periodic differential
operators. Approximation of solutions in the Sobolev class  $H^1(\mathbb{R}^d)$},
Algebra i Analiz  {\bf 18} (2006), no.~6, 1--130; English transl.,
St.~Petersburg  Math. J. {\bf 18} (2007), no.~6, 857--955.

\bibitem[Gr1]{Gr1} Griso~G., \textit{Error estimate and unfolding for periodic homogenization},
Asymptot. Anal. \textbf{40} (2004), 269--286.

\bibitem[Gr2]{Gr2} Griso~G., \textit{Interior error estimate for periodic homogenization},
C. R. Math. Acad. Sci. Paris \textbf{340} (2005), 251--254.

\bibitem[Zh1]{Zh1}Zhikov V.~V.,
\emph{On the operator estimates in the homogenization theory},
Dokl. Ros. Akad. Nauk \textbf{403} (2005), no. 3, 305-308;
English transl., Dokl. Math. \textbf{72} (2005), 535--538.


\bibitem[Zh2]{Zh2}Zhikov V.~V., \textit{On some estimates of homogenization theory},
Dokl. Ros. Akad. Nauk \textbf{406} (2006), no. 5, 597-601; English transl., Dokl. Math. \textbf{73}
(2006), 96--99.


\bibitem[ZhKO]{ZhKO}Zhikov V.~V., Kozlov S.~M., Olejnik O.~A., \textit{Homogenization of
differential operators}, "Nauka", Moscow, 1993; English transl., Springer-Verlag, Berlin, 1994.

\bibitem[ZhPas]{ZhPas} Zhikov~V.~V., Pastukhova~S.~E.,
\textit{On operator estimates for some problems in homogenization theory},
Russ. J. Math. Phys. \textbf{12} (2005), no. 4, 515-524.

\bibitem[McL]{McL} McLean W., \textit{Strongly elliptic systems and boundary integral equations},
Cambridge: Cambridge Univ. Press, 2000.

\bibitem[Pas]{Pas} Pastukhova~S.~E., \textit{On some estimates in homogenization problems
of elasticity theory}, Dokl. Ros. Akad. Nauk \textbf{406} (2006), no. 5, 604-608;
English transl., Dokl. Math. \textbf{73} (2006), 102--106.


\bibitem[Su]{Su} Suslina~T.~A., \textit{Homogenization of the elliptic Dirichlet problem}:
\textit{operator error estimates in} $L_2$, Preprint.

\end{thebibliography}
\end{document}